\documentclass{amsart}
\usepackage{amsfonts,amscd}

\newtheorem{theorem}{Theorem}[section]
\newtheorem{lemma}[theorem]{Lemma}
\newtheorem{corollary}[theorem]{Corollary}
\newtheorem{proposition}[theorem]{Proposition}
\theoremstyle{remark}

\theoremstyle{definition}
\newtheorem{definition}[theorem]{Definition}

\numberwithin{equation}{section} \makeatother

\DeclareMathOperator{\Cdb}{{\mathbb C}}

\DeclareMathOperator{\Ndb}{{\mathbb N}}

\begin{document}
 \textwidth=13.5cm
\textheight=23cm
 \hoffset=-1cm
  \baselineskip=17pt
\title[Open partial isometries and operator positivity]{Open partial isometries and
positivity in operator spaces}

\author{David P. Blecher}
\address{Department of Mathematics, University of Houston, Houston, TX
77204-3008}
\email[David P. Blecher]{dblecher@math.uh.edu}
 \author{Matthew Neal}
\address{Department of Mathematics,
Denison University, Granville, OH 43023}
\email{nealm@denison.edu}
\thanks{*Blecher was partially supported by grant DMS 0400731 from
the National Science Foundation. Neal was supported
by Denison University.}


\begin{abstract}
We first study positivity in C*-modules using tripotents (= partial
isometries) which are what we call {\em open}. This is then used to
study ordered operator spaces via an `ordered noncommutative Shilov
boundary' which we introduce.
  This boundary satisfies the usual universal diagram/property
of the noncommutative Shilov boundary, but with all the `arrows'
completely positive. Because of their independent interest, we also
systematically study open tripotents and their properties.
\end{abstract}

\maketitle

\section{Introduction}

We are interested here in cones of positive
operators $X_+ = \{ x \in X : x \geq 0 \}$, for a space $X$
of bounded linear operators on a Hilbert space,
where $\geq$ denotes the usual order of such operators.
Besides the
intrinsic interest of such objects (for example, operator positivity
plays a central role in many areas of mathematical physics today),
our work is a sequel to \cite{BW}, which was a first step in a new approach
to positivity in an operator space $X$, namely studying it in terms of the
`noncommutative Shilov boundary' of $X$ (see \cite{Ar,Ham,BLM}).
 The latter object is
a Hilbert C*-module, or, equivalently, a {\em ternary ring of
operators} (or {\em TRO} for short), by which we will mean a closed
subspace $Z$ of a C*-algebra $A$ such that $Z Z^* Z \subset Z$. If
$X$ contains positive operators, then so will any containing TRO.
The starting point of the present investigation, was the question of
whether, in this case, all morphisms in the universal property of
the noncommutative Shilov boundary can also be chosen to be positive
(allowing this boundary to be used as a new tool in the study of
ordered operator spaces)? To answer this, one is led immediately to
study positivity in TROs, and we address this topic first (the last
section of our paper concerns positivity in general operator
spaces). In \cite{BW} we considered the case of selfadjoint TROs $Z$
in C*-algebra $A$. In the first part of the present paper, we are
able to generalize, to arbitrary TROs, a fundamental correspondence
from \cite{BW}: we show that the {\em natural cones} in a TRO,
namely $Z \cap A_+$ in the notation above, are in a bijective
correspondence with tripotents (= partial isometries) which are {\em
open}\footnote{These are not the same as the open partial isometries
of \cite{ER3}.} in the sense of \cite{BW}. The emphasis we place on
the relation between positivity and the underlying algebra has its
philosophical origin in \cite{Ef}.  Open tripotents generalize the
notion of `open projections' in C*-algebra theory \cite{P},  which
in turn generalize the notion of `open sets' in topology.    Since
there appears to be no theory of general open tripotents (in our
sense)
 in the literature, we give a careful development of this topic here.
We also briefly discuss {\em compact tripotents}, a notion
which has been treated in the literature in a more general setting
(see e.g.\ \cite{AkeP,BFMP,ER3,ER4,FP}).
 We believe that these objects
should play a role in operator space theory in the future,
in view of the importance of TROs in
that subject (see
e.g.\ \cite[Chapter  8]{BLM} and references therein).
For example, it has strong relations with the
recent study of peak projections and peak tripotents \cite{H,BHN,Betal}.
  In any case, our paper, like its predecessor, in some sense
`marries' the notion of `positivity of Hilbert space operators' to
ideas from the basic structure  theory of JBW*-triples.

Section 4 is mostly devoted to  maximal orderings on TROs. For
example, we analyze a conjectured characterization of maximal
operator space orderings on $*$-TROs from \cite{BW}. Indeed, 1)\ we
show that the proposed characterization is not true for all
$*$-TROs, and 2)\  we isolate the precise class of  $*$-TROs for
which the conjectured characterization is true in general (we call
these the {\em completely orderable} $*$-TROs).

In Section 5, we apply some of our theory from earlier sections to
construct, for an ordered operator space $X$, an ordered version of
the noncommutative Shilov boundary of $X$.  More particularly, we
assign to the usual noncommutative Shilov boundary of $X$, the
natural cone associated with an open tripotent, which in turn is a
supremum of certain `range tripotents' studied in Section 3. This
`ordered boundary' answers the question raised at the start of this
paper: it satisfies the usual universal diagram/property of the
noncommutative Shilov boundary of $X$, but all the `arrows' are
completely positive. We usually do not assume, unlike in the
predecessor \cite{BW} and in the companion paper \cite{BNW}, that
$X$ has an involution $*$. This is simply because of the greater
generality and freedom available in our framework, and because the
analoguous results in the involutive case  are in some sense just a
special case  (with some exceptions that are discussed in
\cite{BNW}). Our results yield, for example, a very algebraic
characterization of the possible `operator space orderings' on a
given operator space (see e.g.\ Theorem \ref{centl}), and
interesting facts about such orderings which are maximal. The
results are particularly good for spaces $X$ whose positive cone
densely spans $X$, which is a common assumption in the theory of
ordered vector spaces.  Indeed it is often  a very reasonable
assumption since `order theory' can say very little about elements
not in the span of the cone.
 In any case, it seems to be
true that ordered
operator spaces with densely spanning cones, which as far as we know have
not hitherto been
considered in the literature,
constitute a setting to which much of the theory of
(unital) operator systems generalizes in a natural and satisfactory
way.  We initiate the study of such spaces
here and in the sequel \cite{BNW}, where, for example, we obtain a new
`unitization' of such spaces (which is universal in that it has the biggest
possible positive cone, as opposed to the unitization from \cite{OSWU,WW,K2}
which has the smallest), and a striking
`rigidity property' (see the end of Section 2 in \cite{BNW}).
This class certainly deserves further study in the future.

We now turn to precise definitions and notation.    Any unexplained
terms below can probably be found in \cite{BLM}, or
any of the other recent books on operator spaces.  All vector spaces
are over the complex field $\Cdb$.
A given cone in a space $X$ will often be written as $X_+$, and we
write $\geq$ for the associated `ordering': $x \geq y$ iff $x - y \in X_+$.
Indeed, we will use the terms `cone' and `ordering' somewhat interchangeably.
A {\em matrix cone} ${\mathfrak c}$ for us will
simply be a sequence $({\mathfrak c}_n)$, where
${\mathfrak c}_n $ is a cone in $M_n(X)$, such that
if $[x_{ij}] \in {\mathfrak c}_n$ then $x_{ii} \in {\mathfrak c}_1$.
A linear map $T : X \to Y$ between spaces with cones
is {\em positive} if $T(X_+) \subset Y_+$.
If the matrix spaces $M_n(X)$ and $M_n(Y)$ also
each have a given cone, for each $n \in \Ndb$,
and if the  canonical `amplification'
$T_n : M_n(X) \to M_n(Y)$ is positive for each $n \in \Ndb$,
then we say that $T$ is {\em completely positive}.
A (resp.\ complete) {\em order embedding} is a (resp.\
 completely) positive map $T$ such that
$T^{-1}$ is (resp.\ completely) positive on Ran$(T)$.
An {\em operator space ordering} or
{\em operator space cone} on an operator space $X$
is a  specified  matrix cone ${\mathfrak c} =
({\mathfrak c}_n)$ so
that there exists a
complete isometry $T$
from $X$ into a C*-algebra $B$, which is at least
completely positive.
That is, $T$ is a complete isometry with
$T_n({\mathfrak c}_n) \subset M_n(B)_+$ for all $n \in \Ndb$.
Of course, it is more natural in some sense to
strengthen this last definition by also
 requiring $T$ to be a complete order embedding;
and we remark that the ordered spaces satisfying this strengthened definition
 have been abstractly characterized
in \cite{OSWU} as the `matrix ordered operator spaces' whose
matrix norms coincide with the `modified numerical radius' norms
(this follows from e.g.\ Corollary 4.11 there).
Nonetheless, the convention
we adopt seems to fit better with our results.  Moreover,
by the characterization from \cite{OSWU} just mentioned, one can
easily see that the two definitions actually coincide
for operator space cones (in our sense) also satisfying the mild
conditions in \cite[Definition 3.3]{OSWU}.

We will sometimes be sloppy, and
use  ${\mathfrak c}$ interchangeably for ${\mathfrak c}_1$ and for the
entire collection $\{ {\mathfrak c}_n \}$, and vice versa.
Similarly, $T({\mathfrak c})$ often denotes
$(T_n({\mathfrak c}_n))$.    We say that one
ordering on $X$ is {\em majorized} by another ordering if the positive cones
for the first ordering are contained in the positive cones for the second
ordering.  We write $X'$ for the dual Banach space
(resp.\ dual operator space) of a Banach space
(resp.\ operator space) $X$, and regard $X \subset X''$.

We refer to e.g.\ \cite{Ham,BLM} for the basic theory of TROs.
TROs were characterized as operator spaces in \cite{NeRu},
but we shall not need this here.
A {\em ternary morphism} on a TRO $Z$ is a linear map $T$ such that
$T(x y^* z) = T(x) T(y)^* T(z)$ for all $x, y, z \in Z$.
A {\em tripotent} is an element $u \in Z$ such that $u u^* u = u$.
We order tripotents by $u \leq v$ if and only if $u v^* u = u$.
This turns out to be equivalent to $u = v u^* u$, or to $u =  u u^* v$,
and implies that $u^* u \leq v^* v$ and  $u u^* \leq v v^*$ \cite{Bat}.
   A {\em WTRO} is a weak*
closed TRO in a W*-algebra. We write $L(Z)$ for the linking
C*-algebra of a TRO, this has `four corners' $Z Z^*$, $Z, Z^*,$ and
$Z^* Z$.   Here $Z Z^*$ is the closure of the linear span of
products $z w^*$ with $z, w \in Z$, and similarly for $Z^* Z$.  If
$E$ is a WTRO, then we write $\bar{L}(E)$ for the W*-algebra linking
algebra, this has `four corners' $\overline{E E^*}^{weak*}$, $E,
E^*,$ and $\overline{E^* E}^{weak*}$. The second dual of a TRO $Z$
is a WTRO, and $L(Z)'' = \bar{L}(Z'')$ (see e.g.\ the proof of
\cite[8.5.17]{BLM}). We will denote by $I$ the injection from $Z$
into $L(Z)$ given by $I(z)=z \otimes e_{12} + z^{\ast} \otimes
e_{21}$. Note that $I(z)^2 = z z^{\ast} \otimes e_{11} + z^{\ast} z
\otimes e_{22}$. For a tripotent $u$, we set  $\hat{u} = \frac{1}{2}
(I(u)+I(u)^{2})$, and  $\breve{u} = \frac{1}{2} (-I(u)+I(u)^{2})$;
these are projections. Define $\Theta : L(Z) \rightarrow L(Z)$ to be
the period 2 $*$-automorphism which changes the sign of the off
diagonal entries. Note that $\breve{u} =
\Theta^{\prime\prime}(\hat{u})$. We say that  a projection $r\in
L(Z)^{\prime\prime}$ is  {\em antisymmetric} if $r \perp
\Theta^{\prime\prime}(r)$, or equivalently, if $r = \hat{v}$ for a
tripotent $v \in Z''$ (see Lemma \ref{adon}).

A subTRO of a TRO $Z$ is a closed subspace of $Z$ which is
closed under the ternary product $x y^* z$.  We write
$\langle {\mathcal S} \rangle$ for the smallest subTRO
containing a given subset ${\mathcal S}$ of $Z$.
An {\em inner ideal} (resp.\ {\em ternary ideal}) of $Z$
is defined to be a closed subspace $J$ with $J Z^* J \subset J$
(resp.\ $J Z^* Z \subset J$ and $Z Z^* J \subset J$).
Clearly inner and ternary ideals are subTROs.
A $*$-TRO is a selfadjoint TRO $Z$ in a C*-algebra $B$,
and by an inner $*$-ideal or ternary $*$-ideal we mean
an inner or ternary ideal which is selfadjoint (that is, closed
under the involution).  A tripotent in a $*$-TRO is
{\em selfadjoint} if $u = u^*$, and {\em central} if $u z = z u$
for all $z \in Z$.

  The {\em Peirce $2$-space} of a tripotent $u$ in
a TRO $Z$ is the subset
$$Z_2(u) = \{ z \in Z :
z = u u^* z u^* u \} =  u u^* Z u^* u  = u Z^* u .$$
Clearly $Z_2(u)$ is an inner ideal of $Z$,
and if  $Z$ is a WTRO then it is weak* closed.
There is a natural product (namely $x \cdot y = x u^* y$) and involution
(namely $x^\sharp =  u x^* u$) on $Z_2(u)$ making the
latter space into a unital $C^*$-algebra.  The identity
element of course is $u$.  If $u \leq v$ then
$Z_2(u)$ is a hereditary C*-subalgebra of $Z_2(v)$, and $u$ becomes
a projection in the last algebra.
If $Z$ is a WTRO then $Z_2(u)$ is a W*-algebra.
The positive cone in the C*-algebra  $Z_2(u)$ will be written as
${\mathfrak c}_{u}$.   Strictly speaking,
we should probably write ${\mathfrak c}^Z_{u}$ for this cone, but to
avoid excessive notation we will write the simpler expression.
We leave it to the reader to make sense of the space
which ${\mathfrak c}_{u}$ lives in (it will always be the
TRO that $u$ belongs to).     It is easy to
check that $u^* Z_2(u)$ is a  C*-subalgebra of $Z^{\ast}Z$, and
the map $z \mapsto u^* z$ is a $*$-isomorphism
from $Z_2(u)$, with the product and involution above,
 onto this C*-subalgebra.  From this it is easy to see that
$${\mathfrak c}_u = \{ z \in Z_2(u) : u^* z \geq 0 \} =
\{ z \in Z : u^* z \geq 0 , z = u z^* u \},$$ and also equals $\{ u
u^* z u^* u z^* u : z \in Z \}$, where these inequalities are in the
C*-algebra $Z^* Z$. If $u \in Z''$, we define ${\mathfrak d}_u$ to
be the cone ${\mathfrak c}_u \cap Z$ in $Z$.    We also write
${\mathfrak c}'_u$ for the weak* closure of ${\mathfrak d}_u$ in
$Z''$.  In contrast to ${\mathfrak c}_u$, the cone ${\mathfrak d}_u$
lies in $Z$, and not in the space $Z''$ which $u$ lies in in this
case. Finally, we will write $Z(u)$ for $Z^{\prime\prime}_{2}(u)
\cap Z$. Following \cite{BW}, we say that a tripotent $u$ in the
WTRO $Z''$ is {\em open}, if when we consider $Z''_2(u)$ as a
W*-algebra in this way, then $u$ is the weak* limit in $Z''$ of an
increasing net from ${\mathfrak d}_u = {\mathfrak c}_{u} \cap Z$.
Beware that this definition differs from the one given in \cite{ER3}
(for example, all unitaries are open in the sense of that paper). We
will show that the spaces  ${\mathfrak d}_{u}$, for open tripotents
$u$, are exactly the natural cones in $Z$, and that this sets up an
order preserving bijection between open tripotents and natural
cones.

\begin{lemma} \label{trp}  A positive  ternary morphism between C*-algebras is
a $*$-homomorphism, and hence it is completely positive.
In particular, a positive linear completely isometric
surjection between C*-algebras is a $*$-isomorphism.
\end{lemma}

\begin{proof}  The first assertion
may be found in the proof of
\cite[Corollary 4.3 (2)]{BW}.
The second we shall not need (it is stated as
background), and it follows from the
well known fact that  the surjective linear
complete isometries between TROs are exactly the ternary isomorphisms.
\end{proof}

\section{Open tripotents and natural cones}

We begin with the following simple but fundamental observation,
which we have not seen in the literature:

\begin{lemma} \label{cruc}  Let $Z$ be a TRO inside a C*-algebra $A$.  Then the
subspace $J(Z) = Z \cap Z^* \cap Z^* Z \cap Z Z^*$ is a
C*-subalgebra of $A$ which is also an inner ideal in $Z$.   Moreover,
the positive cone $J(Z)_+$ of this C*-subalgebra equals $Z \cap A_+$.
  \end{lemma}

\begin{proof}  The  proof is left to the reader, except for the
last assertion.  Clearly $J(Z)_+ \subset Z \cap A_+$.  Conversely,
if $x \in Z \cap A_+$ then of course $x \in Z^*$.
Also, $x^2 \in Z^* Z$, so that $x \in Z^* Z$ since square roots remain in
a C*-algebra.  Similarly, $x \in Z Z^*$, so that $x \in J(Z)$.
\end{proof}

The positive cone $J(Z)_+$ will be called a {\em natural cone}  for
$Z$, and the corresponding ordering on $Z$ is called a {\em natural
ordering}. Since $M_n(Z)$ is a TRO in $M_n(A)$, of course one has a
sequence of cones, $M_n(J(Z))_+ = J(M_n(Z))_+$, but since the cone
$J(Z)_+$ determines the others (see Corollary \ref{acp}), it will
not often be necessary to mention these other cones.  Thus we often
suppress the obvious facts concerning them (in  Section 4 we will
start to be more careful in this regard). We also use the term
`natural  cone' even when the C*-algebra $A$ is not in evidence.
Thus,
 a cone ${\mathfrak d}$ in $Z$ is natural if there
exists a one-to-one ternary  morphism $\varphi : Z \to B$, for a
C*-algebra $B$, such that $\varphi({\mathfrak d}) = \varphi(Z) \cap
B_+$. {\em  Natural dual cones}  for a WTRO $E$ are defined
analoguously (a weak* closed cone such that there exists a
one-to-one weak* continuous ternary  morphism into a W*-algebra
satisfying $\varphi({\mathfrak d}) = \varphi(Z) \cap B_+$. If $Z$ is
a WTRO in a W*-algebra $M$, then $J(Z) = Z \cap Z^* \cap
\overline{Z^* Z}^{weak*}  \cap \overline{Z Z^*}^{weak*}$, a
W*-subalgebra of $M$. To see this, note that the latter space is a
W*-subalgebra, and if $x$ is positive there then $x = (x^*
x)^\frac{1}{2} \in Z^* Z$. Similarly, $x \in Z Z^*$, and so $x \in
J(Z)$.

\medskip

{\bf Remark.}  There seems to be no way to reduce
the study of cones on TROs to the $*$-TRO case studied in
\cite{BW}.  Clearly if $Z$ is a TRO in a C*-algebra
$A$, then $W = Z \cap Z^*$ is a $*$-TRO, and
$W \cap W^2 = J(Z)$ and $W \cap A_+ = Z \cap A_+$.  However the
space $Z \cap Z^*$ depends crucially on the particular ambient
C*-algebra $A$ chosen.
That is, if $\theta  : Z \to B$ is a ternary
morphism and complete order embedding into another C*-algebra $B$, there is no
nice relation in general between $W$ and $\theta(Z) \cap \theta(Z)^*$.

\begin{corollary} \label{acp}  If $\theta : Z \to W$
is a ternary morphism between subTROs of C*-algebras, and if
$\theta$ is positive, then $\theta$ is completely positive.
\end{corollary}

\begin{proof}   Clearly $\theta$ is positive iff
$\theta_{\vert J(Z)}$ is positive as a map from $J(Z)$ to $J(W)$. By
Lemma \ref{trp}, $\theta_{\vert J(Z)}$ is completely positive. Thus,
$\theta_n$ is  positive  as a map from $J(M_n(Z)) = M_n(J(Z))$ to
$J(M_n(W)) = M_n(J(W))$.
\end{proof}

\begin{proposition} \label{3}  Let $Z$ be a TRO, and $u$ a tripotent in
$Z''$.  Then ${\mathfrak d}_u$ is a natural cone in $Z$.
\end{proposition}

\begin{proof}    We define two one-to-one ternary morphisms $\theta$ and $\pi$ from
$Z$ into
$L(Z'')$
as follows:
$$\theta(z) =  \left[ \begin{array}{ccl}
0 & (1-u u^*) z \\ 0 & u^* z \end{array} \right]  \; , \; \; \; \;
\; \; \pi(z) =  \left [ \begin{array}{ccl} z u^* & z (1- u^* u)  \\
0 & 0 \end{array} \right]  \; , \qquad z \in Z  . $$ Let $T(z) =
\theta(z) \oplus \pi(z) \in L(Z'') \oplus^\infty L(Z'')$, this is a
one-to-one ternary morphism.  If $T(z) \geq 0$ then clearly $(1-u
u^*) z = z (1- u^* u) = 0$ so that $z \in Z''_2(u)$. Since $u^* z
\geq 0$ we have $z \in Z \cap {\mathfrak c}_u = {\mathfrak d}_u$.
Conversely, if $z \in {\mathfrak d}_u$ it is even easier to see that
$T(z) \geq 0$.  Thus we have a one-to-one ternary morphism into a
C*-algebra which is an order embedding.
 \end{proof}

\begin{lemma} \label{x}  If $u$ is a tripotent in a TRO $Z$,
and $x \in {\mathfrak c}_u$, then $u^* x = |x| \in Z^* Z$.
\end{lemma}

\begin{proof}
We have  $u^* x \geq 0$ and $u^* x u^* x = (x^\sharp)^* x = x^* x$.  The
result then follows from the uniqueness of square roots.
\end{proof}

We recall that ${\mathfrak c}_u'$ is the weak* closure of ${\mathfrak d}_u
= {\mathfrak c}_u \cap Z$ in $Z''$.

\begin{lemma} \label{supppf}  For every open tripotent
$u \in Z^{\prime\prime}$, we have that
$Z(u)$ is a
weak* dense C*-subalgebra of
$Z^{\prime\prime}_{2}(u)$, and is an inner ideal of $Z$.
Also, ${\mathfrak c}_u' = {\mathfrak c}_u$.
Conversely, every inner ideal of $Z$ which is ternary isomorphic
(or equivalently, completely isometrically isomorphic) to a
C*-algebra arises this way.
\end{lemma}

\begin{proof}
Let $A = Z(u)
= \{ z \in Z : z = u u^* z u^* u \}$, which is an inner ideal.
Let $(u_t)$ be a positive
net in  $A$
converging weak* to $u$.
As remarked earlier $z \mapsto z u^*$ is a $*$-homomorphism on
$Z^{\prime\prime}_{2}(u)$,
and so $u_t  u^* \geq 0$.   By a variant of Lemma \ref{x},
it follows that $u_t  u^*  \in Z Z^*$.  If $z \in A$
then $u_t u^* z \to u u^* z = z$ weak*.   On the
other hand, $u_t u^* \in Z Z^*$ as mentioned above,
so that $u_t u^* z \to z$ weakly in $Z$.  Thus
convex combinations of $u_t u^* z$ converge to $z$ in
norm.  It follows that convex combinations of $u^* u_t u^* z$
converge to $u^* z$ in norm,
and $u^* u_t u^* z
= (u_t^\sharp)^* z = u_t^* z \in Z^* Z$,
so that $u^* z \in Z^* Z$.  Thus
$A$ is a subalgebra
of $Z^{\prime\prime}_{2}(u)$.
To see that it is a $*$-subalgebra, note that
a similar argument to the above shows that $u z^* \in
Z Z^*$.  Since $z^*$ is a norm limit of
convex combinations of $z^* u u_t^*$,
we see that $u z^* u$ is a norm limit of
convex combinations of $u z^* u u_t^* u
= u z^* u_t$.  The latter is in $Z$ since $u z^* \in
Z Z^*$.

Suppose that $\eta \in Z_2''(u)$, and  $z_\lambda \to \eta$ weak*,
with $z_\lambda \in Z$. Then $u_s z_\lambda^* u_t \in A$, since $A$
is an inner ideal. It follows that $u z_\lambda^* u =
z_\lambda^{\sharp}$ is in the weak* closure of $A$.  Hence also
$\eta^{\sharp} \in A^{\perp \perp}$.  Thus $A$ is weak* dense in
$Z_2''(u)$, and so $A'' = Z_2''(u)$ as  von Neumann algebras. Thus
if $\eta \in {\mathfrak c}_u$ then there is an increasing net in
${\mathfrak d}_u$ with weak* limit $\eta$.  So ${\mathfrak c}_u' =
{\mathfrak c}_u$.

Conversely,
assume that $I$ is an inner ideal of $Z$ which is ternary isomorphic
to a C*-algebra
$A$ via a ternary isomorphism
$\psi:A \to I$.  Then $\psi'' : A'' \to I''
= I^{\perp \perp}$ is a one-to-one
ternary morphism.  If $u = \psi''(1)$ then $u$ is a `unitary' tripotent
in $I$, in the sense that $I^{\perp \perp}=I^{\perp \perp}_{2}(u)$, and it is easy to see that
$I^{\perp \perp} = Z^{\prime\prime}_{2}(u)$
since $I^{\perp \perp}$ is an inner
ideal of $Z''$.  Thus $I = Z(u)$.
Moreover, it is clear that $\psi''$ is a
$*$-isomorphism with respect to the canonical product on
$Z^{\prime\prime}_{2}(u)$.  Since the identity of $A''$ is open,
 it is evident that $u$ is open.
\end{proof}

We define an {\em inner C*-ideal} of a TRO $Z$ to be
an inner ideal $J$ of $Z$ with a specified positive cone $J_+$,
which is ternary isomorphic
to a C*-algebra via an order isomorphism.
If, further, $J$ is weak* dense in $Z^{\prime\prime}_{2}(u)$
for a tripotent
$u \in Z''$ such that $J_+ \subset {\mathfrak c}_u$,
then we say that $u$ is a {\em support tripotent} for $J$.

\begin{lemma} \label{supf2}  An inner C*-ideal $J$ of a TRO
has a unique support tripotent $u$, which is automatically open,
and $J = Z(u)$ and
$J_+ = {\mathfrak d}_u$.
\end{lemma}

\begin{proof}
That there exists a support tripotent $u$, which is open, is
proved in the last Lemma.
The proof shows that $J = Z(u)$, and $J_+ = {\mathfrak d}_u$.
 For the uniqueness of $u$, note that if $v$ were
another support tripotent of $J$, then ${\mathfrak c}_u
= {\mathfrak c}_u' \subset {\mathfrak c}_v$,
since $J_+ = {\mathfrak d}_u \subset {\mathfrak c}_v$
and the latter set is weak* closed.
By the well known equivalence of (i) and (ii) in the
next proposition, $u \leq v$.  On the other hand,
$v \in Z''_2(v) = J^{\perp \perp} = Z''_2(u)$, so that
$v = v u^* u = u$.
 \end{proof}

It follows immediately from the last results and discussion earlier
in this section, that every natural cone
on a TRO $Z$ gives rise to an open tripotent, namely the
support tripotent of $J(Z)$.
Conversely, every open tripotent gives a natural cone
by Proposition \ref{3}.

For the following result, we recall the definition
$\hat{u} = \frac{1}{2} (I(u)+I(u)^{2})$ from the Introduction.

\begin{proposition} \label{prop3}  For tripotents $u, v$ in a TRO $Z$,
 the following are equivalent:
\begin{itemize}
\item [(i)] $u \leq v$ in $Z$.
\item [(ii)]  ${\mathfrak c}_{u} \subset {\mathfrak c}_{v}$.
  \item [(iii)]  $Z_2(u)$ is a C*-subalgebra of $ Z_2(v)$.
\item [(iv)]  $u \in Z_2(v)$ and $u$ is a projection in
that C*-algebra.
\item [(v)]   $\hat{u} \leq \hat{v}$.
\end{itemize}
\end{proposition}

\begin{proof}  These are all essentially well known (see
e.g.\ \cite{Bat}), and
easy exercises, except perhaps the equivalence with (v).
If $\hat{u} \leq \hat{v}$, then $uu^{\ast} \leq vv^{\ast}$ and
$u^{\ast}u \leq v^{\ast}v$.  Since $\hat{u}\hat{v}=\hat{u}$, we see that
$u/2=uu^{\ast}v/4+uv^{\ast}v/4$ and, thus, $u=uu^{\ast}v$ and $u \leq
v$. The other direction is obtained by multiplying $\hat{u}$ by
$\hat{v}$ and using the equations $uu^{\ast}v=vu^{\ast}u=u$.
\end{proof}

\begin{corollary}  \label{tprop3}  If $u, v$ are open
tripotents in the second dual of a TRO $Z$, then the following are equivalent:
\begin{itemize}
\item [(i)] $u \leq v$ in $Z''$.
\item [(ii)]  ${\mathfrak d}_{u} \subset {\mathfrak d}_{v}$.
\item [(iii)]  $Z(u)$ is a C*-subalgebra of $Z(v)$.
\end{itemize}
Also, the correspondence established above
between natural cones on $Z$ and open tripotents
in $Z''$ is bijective, and preserves `order' (ordering cones by inclusion).
 \end{corollary}

\begin{proof}  Suppose that ${\mathfrak d}_{u} \subset {\mathfrak d}_{v}$.
Taking weak* closures, ${\mathfrak c}_{u} \subset {\mathfrak c}_{v}$,
so that $u \leq v$ by Proposition \ref{prop3}.
We leave the other implications as an exercise, using of course
that proposition and the
earlier established facts summarized in the remark above Proposition \ref{prop3}.
  \end{proof}

{\bf Remarks.}   1) \ Variants of the arguments above
show that there is a bijective order preserving correspondence between natural
dual cones in a WTRO $E$, and tripotents in $E$.

2) \ It follows from Corollary \ref{tprop3}
that maximal natural cones for $Z$
will correspond to maximal open tripotents in $Z''$.
Maximal open tripotents in $Z''$ are studied in
the second half of Section 4 below.

\begin{lemma} \label{uxinE}   Let $Z$ be a
 TRO, and let $E = Z''$, also a TRO in the canonical way.
Let $u$ be a tripotent  in $E$, so that ${\mathfrak d}_u$  is a natural cone
by Proposition {\rm \ref{3}}.
We have:
\begin{itemize} \item [(1)]  ${\mathfrak c}'_u \subset {\mathfrak c}_u$.
\item [(2)]   ${\mathfrak c}'_u = {\mathfrak c}_v$ for an
open tripotent $v \in E$ with $v \leq u$.
\item [(3)]  The closed span $J_u$ of
${\mathfrak d}_u$  is a C*-subalgebra of $E_2(u)$, the latter
regarded as a C*-algebra in the canonical way.
Also, ${\mathfrak d}_u$  is
the positive cone of this C*-algebra $J_u$.
\end{itemize}
\end{lemma}

\begin{proof}   Item (1) is obvious, and
(2) follows from
Proposition \ref{prop3}, and Lemmas \ref{supppf} and
\ref{supf2}.   For  (3) note that
${\mathfrak d}_u = {\mathfrak d}_v$,
and so $J_u$  is a C*-subalgebra of $E_2(v)$,
which in turn is a C*-subalgebra of $E_2(u)$
by Proposition \ref{prop3} (iii).
\end{proof}

{\bf Remark.}   We do not know if it is true that if
 $Z$ is a TRO in a C*-algebra $A$,
and if $Z''$ is regarded as a  TRO in $A''$, then
 $J(Z)$ is weak* dense in $J(Z'')$.

\begin{theorem} \label{chopt}  Let $Z$ be a TRO, set $E = Z''$, and
let $u$ be a tripotent in $E$.
The following are equivalent: \begin{itemize}
\item [(i)]  $u$ is an open tripotent (i.e.\ there is a net $(x_t)$
in $Z$ converging weak* to $u$,
satisfying: $u^* x_t \geq 0$, $u x_t^* u = x_t$ for all $t$, and
$(u^* x_t)$ is an increasing net).
\item [(ii)]  $u \in {\mathfrak c}'_u$.
\item [(iii)]  $\hat{u} = \frac{1}{2} (I(u)+I(u)^{2})$ is an open
projection in $L(Z)^{\prime \prime}$.
\item [(iv)]  $\breve{u} = \frac{1}{2} (-I(u)+I(u)^{2})$ is an open
projection in $L(Z)^{\prime \prime}$.
\item [(v)]  $u$ is a support tripotent for an inner C*-ideal in $Z$.
 \item [(vi)]  ${\mathfrak c}'_u = {\mathfrak c}_u$.
\item [(vii)]  The closed span $J_u$ of
${\mathfrak d}_u$  is weak* dense in $E_2(u)$.
 \item [(viii)]    $-u$ is an open tripotent.
\end{itemize}  \end{theorem}

\begin{proof}   Lemmas \ref{supppf} and \ref{supf2}
give (i) $\Leftrightarrow$ (v), and the fact that
(i) implies (vi) and (vii).  It is easy to see the equivalence
of (viii) with (i) from the definition.
Clearly (i) implies (ii), and (vi) implies (ii).

(vii) $\Rightarrow $ (ii) \   By Lemma \ref{uxinE},
$J_u$ is a C*-subalgebra of $E_2(u)$.   If also
$J_u^{\perp \perp} =  E_2(u)$ it follows that
$J_u = E_2(u) \cap Z$, and $E_2(u)$ is the `second dual C*-algebra' of
$J_u$.   Thus $(J_u)_+ = {\mathfrak d}_u$, and
$u \in {\mathfrak c}_u'$.

(ii) $\Rightarrow $ (iii) \ Suppose that
$(x_t)$ is a net  in ${\mathfrak d}_u$ converging w* to $u$.  Let
$$r_t  = \frac{1}{2} \left [ \begin{array}{ccl}
x_t u^* & x_t \\ x_t^* & u^* x_t \end{array} \right] .$$
Since $x_t$ is selfadjoint in $Z^{\prime\prime}_{2}(u)$, we have
$u^* x_t = x_t^* u$ and $x_t u^* = u x_t^*$.  Thus $r_t$ is self-adjoint,
and it clearly converges weak* to $\hat{u}$.   In fact $r_t \in L(Z)$.
Indeed,
$x_t u^* \in Z Z^*$ and $u^* x_t \in Z^* Z$ as in Lemma  \ref{x}.
Then $\hat{u}$ is open, since $r_t \hat{u} = r_t$.

(iii) $\Rightarrow $ (i) \ Suppose that $\hat{u}$ is open.  Let
$$r_t  = \frac{1}{2} \left [ \begin{array}{ccl}
a_t & x_t \\ x_t^* & b_t \end{array} \right] $$
be a positive net in $L(Z)$ increasing up to $\hat{u}$.
We have $\frac{1}{4} a_t u
+  \frac{1}{4} x_t u^* u = \frac{1}{2} x_t$, which,
multiplying by $u^{\ast}u$, yields $x_t = a_t u$, and $x_t = x_t u^* u$.
Similarly $x_t = u u^* x_t$, so that $x_t \in
Z(u)$.    Since $(a_t)$ is positive
and increasing, and $x_t = a_t u$, we have that $(x_t)$ is
positive and increasing
in $Z^{\prime\prime}_{2}(u)$.  Thus $u$ is an open tripotent.

(iii) $\Leftrightarrow $ (iv) \  Follows from the equivalence of
(viii) with (iii).
 \end{proof}

{\bf Remark.} One may prove directly that (ii) implies (i). Indeed
the proof of Lemma \ref{supppf} shows that $A = Z(u)$ is a
C*-subalgebra of $Z_2''(u)$, and $A'' = Z_2''(u)$ as  von Neumann
algebras. It follows that there is an {\em increasing} positive net
in $A$ with weak* limit $u$.  This gives (i).

\begin{corollary}  \label{tprop3b}  If $Z$ is a TRO and $u, v$ are tripotents in
$Z''$, with $v$ open and $u \leq v$, then $u$ is open iff it is open
as a projection in $Z''_2(v)$.
\end{corollary}

\begin{proof}  If
 $z_t \to u$ weak*, with $(z_t)$ an increasing net in
 ${\mathfrak d}_{u}$,
then $z_t$ is an increasing net in ${\mathfrak c}_{v}$.  So
$u$ is an open  projection in $Z''_2(v)$.  Conversely, if
$u$ is an open  projection in $Z''_2(v)$, then
there is a net $x_t \in  {\mathfrak d}_{v}$ with $x_t \leq u$ in $Z''_2(v)$,
and $x_t \to u$ weak*.  We have $u v^* x_t = x_t = u u^* x_t$, and
similarly $x_t u^* u  = x_t$.  Thus $x_t \in Z(u)$, and
indeed $x_t \in {\mathfrak d}_{u}$ since $u^* x_t = v^* x_t \geq 0$.
 Thus $u$ is open  by Theorem \ref{chopt} (ii).
\end{proof}

{\bf Remark.}    The open tripotents in the second dual of a
C*-algebra, which are projections, are exactly the usual open
projections \cite{P}.

More generally, if $u$ is an open tripotent in $Z''$, then it is
easy to see from the proof of (iii) $\Rightarrow$ (i) in Theorem
\ref{chopt}, that $u u^*$ and $u^* u$  are open projections in $(Z
Z^*)''$ and  $(Z^* Z)''$ respectively (cf.\ \cite{ER3}).

 \begin{proposition} \label{incr}  An increasing net $(u_t)$ of
open tripotents in $Z''$
has a least upper bound tripotent $u$ in $Z''$, namely its weak* limit,
and $u$ is also open.  In terms of
cones, ${\mathfrak d}_u$ is the norm closure of the
union of the cones  ${\mathfrak d}_{u_t}$.  That is,
the norm closure of a union of a nested collection of
natural cones is a natural cone.
  \end{proposition}

\begin{proof}   It is well known \cite[Proposition 3.8]{Bat},
and easy to argue directly by
a weak* limit argument using separate weak* continuity of
the product, that the net has an upper bound tripotent $w$.
 Working in $Z''_2(w)$, the $u_t$ become an increasing
net of projections.  Hence
they have a supremum projection $u$, which they converge to strongly.
  It is easy to check that
$u$ is also the supremum as a tripotent in $Z''$.
The spaces $Z''_2(u_t)$ are W*-subalgebras of $Z''_2(u)$.  To see
that $u$ is open, by Theorem \ref{chopt} (ii) it
suffices to show that $u_t$ is in the
weak* closure of ${\mathfrak d}_{u}$.
However $u_t$ is in the
weak* closure of ${\mathfrak d}_{u_t}$, and
$${\mathfrak d}_{u_t} = {\mathfrak c}_{u_t}
\cap Z \subset {\mathfrak c}_{u} \cap Z  = {\mathfrak d}_{u} .$$

Since the $u_t$ are open projections in the C*-algebra
$Z''_2(u)$ by Corollary \ref{tprop3b},
the final assertion is essentially well known.
We include a proof for completeness.  It is clear that
${\mathfrak d}_u$ contains the norm closure of the
union of the cones  ${\mathfrak d}_{u_t}$.  To get the
reverse inclusion, let $A_t$ be the C*-subalgebra
$Z(u_t)$ of $A = Z(u)$.
In turn $A$ is a C*-subalgebra of $Z''_2(u)$,
and indeed $Z''_2(u)$ is the second dual
C*-algebra of $A$.  The positive cones of $A_t$ and $A$ are
${\mathfrak d}_{u_t}$ and ${\mathfrak d}_{u}$
respectively.   The weak* closure of $\cup_t \, Z''_2(u_t)$
equals $Z''_2(u)$, since any $\eta \in Z''_2(u)$
is the strong limit of $u_t u^* \eta u^* u_t$,
and we have $u_t u^* \eta u^* u_t \in  Z''_2(u_t)$.  It follows by basic
functional analysis that the norm closure of  $\cup_t \, A_t$ is $A$.  From this it is
clear that the closure of $\cup_t \, (A_t)_+$ is $A_+$.
\end{proof}

We leave the following as an exercise.

\begin{corollary} \label{dirs}  The `$L^\infty$-direct
sum' $\oplus^\infty_i \; Z_i$ of naturally ordered TROs (resp. dual
naturally ordered WTROs),
 with the
obvious cone,  is again a naturally ordered TRO
(resp.\
dual naturally ordered WTRO).
\end{corollary}

\section{Further properties of open tripotents}

We begin with some facts and lemmas on range tripotents,
almost all of which are well known: in the literature
 (see e.g.\ \cite{AkeP,Bat} and the cited papers of
Edwards and R\"uttimann, especially \cite{ER4}) or folklore.
However, since the arguments are short
and simple we include them here for the readers convenience.
Later in this section we  establish the  basic `calculus'
of open tripotents, following (and freely using ideas from)
the  basic `calculus' of open projections,
and the `calculus' of tripotents established in the
aforementioned papers.
In view of the bijective correspondence
from Corollary \ref{tprop3}, this `calculus' may be reread
as constituting most of the basic `theory of natural cones'.  We will
not usually explicitly state the `cone version' of each result
below, but leave this to the reader.

Let $E$ be a WTRO.
For each $x \in E$, we consider the range tripotent
$r(x)$ in $E$.
This is the tripotent in $E$ with the
property that $x = r(x) |x|$ and $r(x)^{\ast}r(x)$ is
the support projection of $|x|$
(namely, $r(x)$ is the partial isometry  in the polar decomposition
of $x$, see e.g.\ 8.5.22 in \cite{BLM}).
Such a tripotent is unique: it is the
smallest tripotent $u$ in $E$ with the
property that $x = u |x|$.
(To see this, note that if $x = u |x|$, then
$u r(x)^* r(x) |x| = u |x| = x = r(x) |x|$, and so $u r(x)^* r(x)  = r(x)$.)
We have that $$x r(x)^* r(x) = r(x) |x| r(x)^* r(x) = r(x) |x| = x , $$
since $1 - r(x)^* r(x)$ is the projection onto Ran$(|x|)^\perp = {\rm Ker}(|x|)$.
Also, $$|x| r(x)^* r(x) = (r(x)^* r(x) |x|)^* = |x|,$$
and so
$$r(x)x^{\ast}r(x) = r(x) |x| r(x)^* r(x) = r(x) |x| = x .$$
This shows that $x \in Z''(r(x))_+$, since $r(x)^* x = |x| \geq 0$.

Define $C_{0}(x)$ to be the norm closure of the
span of odd polynomials of $x$.
By \cite[Lemma 3.2]{ER3}, $C_{0}(x)$
is a commutative sub-C*-algebra of
$Z''_{2}(r(x))$, and odd polynomials in $x$
 are the same in either product.
 It follows easily that $C_{0}(x)$
is ternary isomorphic to $C_{0}({\rm Sp}(x))$ where the spectrum is taken in
$Z''_{2}(r(x))$. Thus the elements in $C_{0}(x)$ coincide with the
usual functional calculus in the W*-algebra
$Z''_{2}(r(x))$.
If $Z$ is a $*$-TRO and if $x \in Z_{sa}$, then
clearly $r(x)$ is selfadjoint.

\begin{lemma} \label{toadd}  If $u$ is a tripotent in a WTRO $Z$, then
${\mathfrak c}_{u} = \{ z \in Z : r(z) \leq u \}$.
If $Z$ is a TRO
and $u \in Z''$,
then ${\mathfrak d}_{u} = \{ z \in Z : r(z) \leq u \}$.
  \end{lemma}

\begin{proof}
If $r(z) \leq u$ then $z \in {\mathfrak c}_{r(z)} \subset {\mathfrak c}_{u}$,
by Proposition \ref{prop3}.
  Conversely,
if $z \in {\mathfrak c}_{u}$, then $u |z| = u u^* z = z$ by
Lemma \ref{x}, so that $r(z) \leq u$.
The final assertion follows immediately from the first one.
\end{proof}

\begin{lemma} \label{rtrip}
If $u$ is a tripotent in a WTRO $E$, and if $x \in {\mathfrak c}_u$, then
$r(x)$ is the support projection of $x$ in the W*-algebra $E_{2}(u)$.
 \end{lemma}

\begin{proof}
Since $x \in {\mathfrak c}_u$
we have $r(x) \leq u$ by Lemma \ref{toadd}, and so $r(x)$ is
a projection in $E_{2}(u)$.   We have $r(x) u^* x = r(x) r(x)^* x = x$.
If $v$ is another projection in $E_{2}(u)$ with $v  u^* x = x$
then $v u^* r(x) |x| =
r(x) |x|$.
This implies that $v u^* r(x)  = r(x)$ so that $v \geq r(x)$.
Thus $r(x)$ is the support projection of $x$ in $E_{2}(u)$.
\end{proof}

\begin{lemma} \label{Matt}  Let $Z$ be a TRO.  The range tripotent $r(x)$
of any $x \in Z$ coincides with the weak* limit of
$x^{1/(2n-1)}$, the `power' taken in the W*-algebra $Z''_2(r(x))$,
and is open. Furthermore, $Z(r(x))$ is the smallest
inner ideal in $Z$ containing $x$.
\end{lemma}

\begin{proof}
We work in the W*-algebra $Z''_{2}(r(x))$. As we said earlier, $x
\geq 0$ in $Z''_2(r(x))$. The first statement then follows from
Lemma \ref{rtrip}, and well known properties of support projections
in a W*-algebra.  As mentioned at the beginning of this section, odd
polynomials of $x$ are the same in either of the two products, and
lie in $Z$.  Since $x^{1/(2n-1)}$ is a norm limit of odd polynomials
of $x$ in $Z''_{2}(r(x))$, we have that $x^{1/(2n-1)}$ lies in $Z$.
Thus $r(x)$ is open by (ii) of Theorem \ref{chopt}.  For the last
statement, it is clear that any weak* closed inner ideal of $Z''$
containing $x$ must contain $r(x)$, and, thus, all of
$Z''_{2}(r(x))$ by definition of inner ideal. Hence, any inner ideal
$I$ of $Z$ must contain $Z(r(x))$ since $\overline{I}^{weak*}$ is an
inner ideal of $Z''$ (see also see Lemma 3.7 of \cite{ER3}).
\end{proof}

In the last proof we showed that
 if $x \in {\rm Ball}(Z)$, then $r(x)$ is an
increasing weak* limit
of powers $x^{1/(2n-1)}$, which in turn are norm limits of odd
polynomials in $x$.   This is also true if $x \in E$ for  a WTRO $E$,
with the weak* limit being in the weak* topology of $E$.
We will use these facts frequently in the sequel,
often silently.

\begin{corollary} \label{incr2}  A tripotent is open iff it is
a weak* limit of an increasing net of range tripotents.
\end{corollary}

\begin{proof}
A tripotent which is a limit of an increasing net of range tripotents
is open by Proposition \ref{incr} and Lemma \ref{Matt}.  Conversely,
suppose that $u$ is open.
Then $u$ is a  weak* limit of an increasing net
$(z_t)$ in $Z(u)$.  By Lemma \ref{rtrip}, $r(z_{t})$  is
the support projection of $z_t$ in the W*-algebra $Z''_{2}(u)$. Thus, the
net of range tripotents
$(r(z_t))$  are increasing, and we have $z_t \leq r(z_t) \leq u$.
It follows that $r(z_t) \to u$ weak*.
\end{proof}

{\bf Remark.}  One may also construct `open spectral tripotents'
as follows.  If $x$ is an element in a TRO $Z$, then
$x$ is positive in the W*-algebra $Z''_2(r(x))$.  If $U$ is any
open set in Sp$(x) \subset [0,\infty)$, then the
spectral projection $\chi_U(x)$ in the W*-algebra $Z''_2(r(x))$
is open, and hence it is an open tripotent in $Z''$.
A special case of course is if we take $U = (0,\infty)$,
then the associated `spectral  open tripotent' $u$ is just $r(x)$.
Indeed this is a  well known formula for a support projection in a W*-algebra.

We now turn to properties of general open tripotents.

\begin{proposition} \label{imop}  If $\theta : Z \to W$
is a ternary morphism between TROs, and if
$u$ is an open  tripotent in
 $Z''$, then $v = \theta''(u)$ is an open  tripotent in $W''$.
Also, $\theta$ restricts to a
$*$-homomorphism $Z(u) \to Z(v)$, and $\theta({\mathfrak d}_u) \subset
{\mathfrak d}_{v}$.   If also $\theta$ is surjective, then
$\theta({\mathfrak d}_u) = {\mathfrak d}_{v}$.
\end{proposition}

\begin{proof}
Suppose that $z_t \in Z$ with $u^* z_t \geq 0, z_t = u z_t^* u,$ and
$z_t \to u$ weak*.  Applying $\theta$ we obtain a net $(\theta(z_t))$
with analoguous properties, so that $v$ is open by Theorem \ref{chopt} (ii).
Clearly $\theta''$ is a $*$-homomorphism $Z_2''(u) \to Z_2''(v)$,
and thus restricts to a
$*$-homomorphism $Z(u) \to Z(v)$.  If also $\theta$ is surjective, then
$\theta''(Z_2''(u))$ is an inner ideal containing $v$, so that
$\theta''(Z_2''(u)) = Z_2''(v)$.  From this it follows that
$\theta(Z(u)) = Z(v)$ and  $\theta({\mathfrak d}_u) = {\mathfrak d}_{v}$.
 \end{proof}

{\bf Remark.}  Thus $\theta$ will be completely positive on $Z$ with
its ordering determined by $u$.

\begin{lemma} \label{adon}  If $r$ is
an antisymmetric projection in $\bar{L}(E)$
for a WTRO $E$, then $r=\hat{v}$
for a tripotent $v \in E$.
\end{lemma}
\begin{proof}
Suppose that $r=a\otimes e_{11}+b \otimes e_{12} + b^{\ast} \otimes
e_{21} +c \otimes e_{22}$. Squaring $r$ yields $a^{2}+ bb^{\ast}=a$.
Antisymmetry yields $a^{2}-bb^{\ast}=0$. Thus $2a^{2}=a$. Hence $2a$
is a projection. Since $(2b)(2b)^{\ast}=4a^{2}=2a$, we see  $2b$ is
a partial isometry. Similarly $(2b)^{\ast}(2b)=2c$. The result
follows.
\end{proof}

We now turn to the supremum $u \vee v$ of two tripotents.
Most of the following result is in \cite{AkeP}, but for convenience
we give quick proofs.

\begin{lemma}\label{sup} Suppose $u,v$ are tripotents in a WTRO $E$.
The following are equivalent:
\begin{itemize}
\item[(i)] $u \vee v$ exists.
 \item[(ii)] $\hat{u} \perp \breve{v}$.
 \item[(iii)] $uv^{\ast}v=uu^{\ast}v$ and $vv^{\ast}u=vu^{\ast}u$.
\item[(iv)]  $\{ u , v \}$ is dominated by a tripotent $w \in E$.
 \end{itemize}
In this case, $\widehat{u \vee v} = r(\hat{u}+\hat{v})
= \hat{u} \vee \hat{v}$
and $u \vee v =r(u+v)$.
 \end{lemma}
\begin{proof}
(i)  $\Rightarrow$ (iv) \ Obvious.

(iv) $\Rightarrow$ (iii) \  Suppose that $u \leq w$ and $v \leq w$.
We have $$uv^{\ast}v=uw^{\ast}wv^{\ast}v=uw^{\ast}v
=uu^{\ast}uw^{\ast}v=uu^{\ast}v .$$  The proof of the other statement is similar.

(iii) $\Rightarrow$ (ii) \ Clear from direct multiplication, noting that
$uu^{\ast}vv^{\ast}=uv^{\ast}$ and
$u^{\ast}v=u^{\ast}uv^{\ast}v$.

(ii) $\Rightarrow$ (i) \
 Multiplication shows that
$(\hat{u}+\hat{v})\Theta(\hat{u}+\hat{v})=0$.   Hence,
$p(\hat{u}+\hat{v})\Theta(q(\hat{u}+\hat{v}))=0$
for any odd polynomials $p$ and $q$. Taking weak* limits,
it follows that $r(\hat{u}+\hat{v}) \Theta(r(\hat{u}+\hat{v})) = 0$,
and so $r(\hat{u}+\hat{v})$ is
an antisymmetric projection, and thus equals $\hat{w}$
for a tripotent $w$ by Lemma \ref{adon}.  As is well known,
$\hat{u} \vee \hat{v}=r(\hat{u}+\hat{v})$,
so that $\hat{u} \vee \hat{v}= \hat{w}$.
By Proposition  \ref{prop3} we have $u \leq w$ and $v \leq w$,
and so $u \vee v \leq w$.  If $w_0 = u \vee v$
then $\widehat{w_0} \geq \hat{u}$ and
$\widehat{w_0} \geq \hat{v}$, so that $\widehat{w_0} \geq
\hat{u} \vee \hat{v}= \hat{w}$.
Proposition  \ref{prop3} gives $w_0 \geq w$, so that
$w = u \vee v$.
Finally,
$r(u+v)=u \vee v$ by \cite[Proposition 3.9 (i)]{Bat}.
\end{proof}

\begin{proposition} \label{prop4}  Let $Z$ be a TRO.
A family $\{ u_i : i \in I \}$ of open tripotents in
$Z''$, which are bounded above by a tripotent in $Z''$,
has a least upper bound amongst the tripotents in
$Z''$, and this is an open tripotent.
\end{proposition}

\begin{proof}
We first show that if $u, v$ are open tripotents in
$Z''$, which are bounded above by a tripotent,
then  the sup tripotent
$w=u \vee v$ (which exists by the previous lemma)
is open.
By Lemma \ref{uxinE} and its proof,
Span$({\mathfrak d}_w) = Z(e)$ for an open tripotent $e \leq w$ with
${\mathfrak d}_w = {\mathfrak d}_e$.
Since ${\mathfrak d}_u \subset {\mathfrak d}_w
= {\mathfrak d}_e$, we have $u \leq e$.
Similarly $v \leq e$, and so $e = w$.

It now follows by induction that the supremum of any finite family of
open tripotents  which are bounded above by a tripotent,
is open.  The result then follows easily from
Proposition \ref{incr}.
\end{proof}

 We say that tripotents $u$ and $v$ {\em commute}
if $v^* u = u^* v$ and $v u^* = u v^*$ (Harris calls this
 {\em $*$-commuting}).  We say that
$u \perp v$ if these quantities are zero.

\begin{corollary} \label{prop1}   Let $Z$ be a TRO,
and let $u, v$ be two commuting
open tripotents satisfying
$v u^* u = v v^* u$.
The supremum $u \vee v$ in the
set of tripotents in $Z''$ exists, is open, and is given by the
formula $u + v - v v^* u$.  In particular, $u \vee v = u + v$ is open
if $u \perp v$.
\end{corollary}

\begin{proof}
Let $u, v$ be as stated.
Then, by Lemma \ref{sup}, $w=u \vee v$ exists, and it is open
by Proposition \ref{prop4}. By commutativity,
$uw^{\ast}v=uu^{\ast}v=vu^{\ast}u=vw^{\ast}u$, so $u$ and $v$
commute as projections in $Z_{2}(w)$. As is well known in
this case,
$w=u\vee v=u+v-uw^{\ast}v =u+v-uu^{\ast}v$. The result follows.
\end{proof}

{\bf Remark.}  We do not know how to describe cones corresponding
to suprema of tripotents, even if they are orthogonal.
For example, the usual cone in $M_2$ seems
not nicely related to ${\mathfrak d}_{E_{11}}$ and ${\mathfrak d}_{E_{22}}$.

\begin{corollary} \label{prop2}  Let $Z$ be a TRO.
The infimum $u \wedge v$ of two commuting
open tripotents $u$ and $v$ in $Z''$ is open, and is given by the
formula $\frac{1}{2} (v v^* u + v u^* u)$.
\end{corollary}

\begin{proof}
Let $u, v$ be two commuting
open tripotents.  Then $\hat{u}$ and $\hat{v}$ are two
commuting
open  projections in $L(Z)''$.  It is well known
that the infimum of two commuting
open projections is open \cite{Ake,Ake2}, and so
$\hat{u} \hat{v}$ is open.   However
$\hat{u} \hat{v} = \hat{w}$, where $w = \frac{1}{2} (v v^* u + v u^* u)$.
Since $w$ is a tripotent, $w$ is open by Theorem \ref{chopt}.
It is easy to check that $w = u \wedge v$ (or see e.g.\ \cite{AkeP}).
\end{proof}

The infimum of a collection of open tripotents amongst all
tripotents in $Z''$ need not be open.  However there is
an infimum amongst the {\em open} tripotents in $Z''$:

\begin{lemma} \label{intcon}  Suppose that $Z$ is a TRO and that
${\mathcal F} = \{ u_\alpha : \alpha \in I \}$ is a collection of open tripotents in $Z''$.
Then there exists an infimum $u$ for
${\mathcal F}$  amongst the open tripotents in $Z''$.
 Also, ${\mathfrak d}_{u} = \cap_{\alpha \in I} \, {\mathfrak d}_{u_\alpha}$.
\end{lemma}

\begin{proof}   Clearly the infimum is the
supremum of the open tripotents $v$ such that $v \leq u_\alpha$ for every $\alpha
 \in I$.   This is open by Proposition \ref{prop4}.
Since $u \leq u_\alpha$ for every $\alpha
 \in I$, we have ${\mathfrak c}_{u} \subset {\mathfrak c}_{u_\alpha}$,
and so ${\mathfrak c}_{u} \subset  \cap_{\alpha \in I} \, {\mathfrak c}_{u_\alpha}$.
Conversely, if $x \in \cap_{\alpha \in I} \, {\mathfrak c}_{u_\alpha}$
then by  Lemma \ref{toadd} we have $r(x) \leq u_\alpha$ for every $
\alpha \in I$, and so $r(x) \leq u$.  By  Lemma \ref{toadd} again,
$x \in {\mathfrak c}_{u}$.  Thus
${\mathfrak c}_{u} = \cap_{\alpha \in I} \, {\mathfrak c}_{u_\alpha}$,
and the result is now obvious.
\end{proof}

The lemma asserts that an intersection of natural cones
is a natural cone.
This is valid at the matrix level too,
and this will play a role later.  That is,
$({\mathfrak d}_u)_n = \cap_\alpha \; ({\mathfrak d}_{u_\alpha})_n$
for all $n \in \Ndb$.  A direct proof of this: if
$\{ {\mathfrak c}_\alpha \}$ is a family of natural
dual cones in a WTRO $E$, and if $T_\alpha : E \to B_\alpha$ is a
ternary morphism which is an order embedding
for the cone ${\mathfrak c}_\alpha$, then the map
$T : Z \to \oplus^\infty_\alpha \, B_\alpha$
taking $z \in Z$ to $\oplus_\alpha \, T_\alpha(z)$
is a ternary morphism which is an order embedding
for the cone ${\mathfrak c} =
\cap_\alpha \, {\mathfrak c}_\alpha$.  Hence it
is a complete order embedding
by Corollary \ref{acp}.  Thus if $[z_{ij}] \in M_n(E)$
then $[z_{ij}] \in {\mathfrak c}_n$
iff $[T(z_{ij})] \geq 0$ iff $[T_\alpha(z_{ij})] \geq 0$
for each $\alpha$.  In turn, this happens
iff $[z_{ij}] \in ({\mathfrak c}_\alpha)_n$
for each $\alpha$, that is, iff
$[z_{ij}] \in \cap_\alpha \; ({\mathfrak c}_\alpha)_n$.
Thus ${\mathfrak c}_n = \cap_\alpha \; ({\mathfrak c}_\alpha)_n$.
In terms of tripotents,
$$(\wedge_\alpha \, u_\alpha) \otimes I_n =
\wedge_\alpha \; (u_\alpha \otimes I_n) .$$
The same argument works for
natural cones in a TRO $Z$, or this can be deduced
from the above by taking $E = Z''$.

\begin{corollary} \label{inmn}
Let $n \in \Ndb$.
 The map $z \to z \otimes I_n$ from a WTRO $E$ into
$M_n(E)$ is a one-to-one ternary morphism that preserves
infima of tripotents, and
suprema of tripotents where they exist.
\end{corollary}

\begin{proof}  The statement about infima is demonstrated above.
Next, if $\{ u_\alpha \}$ is a family of tripotents in
$E$ which are bounded above by a tripotent $u$,
then $u_\alpha \otimes I_n \leq u  \otimes I_n$,
so that $\{ u_\alpha \otimes I_n \}$ is bounded above
by $u \otimes I_n$,
where $u = \vee_\alpha \; u_\alpha$.
Conversely, suppose that $\{ u_\alpha \otimes I_n \}$ is bounded above
by a tripotent $w \in M_n(E)$.
It is easy to see that $w_{ij} u_\alpha^* u_\alpha
 = u_\alpha$ for each $i, j, \alpha$.  Thus
$u(w_{ij}) \geq u_\alpha$.   Hence $\{ u_\alpha \}$ is bounded above
by $v = \wedge_{i,j} \;  u(w_{ij})$,
so that $u \leq v$ where $u = \vee_\alpha \; u_\alpha$.  Note that
$$w_{ij} u^* u = w_{ij} u(w_{ij})^* u(w_{ij}) u^* u =
u(w_{ij}) u^* u = u .$$
Thus $u_\alpha \otimes I_n \leq u  \otimes I_n \leq w$.
Hence the supremum of $\{ u_\alpha \otimes I_n \}$
amongst the tripotents in $M_n(E)$ is
$u \otimes I_n$.
\end{proof}

We will need a few `matrix tricks' which we have not seen in
the literature.
For $x \in {\rm Ball}(Z)$ define $$\hat{x} = \frac{1}{2}  \left[ \begin{array}{ccl}
|x^*| & x \\ x^* & |x| \end{array} \right]  \in L(Z) .$$
Writing this as a sum of a diagonal matrix and an off-diagonal matrix,
we see that $||\hat{x}|| \leq 1$.   Letting
$$z =  \frac{1}{\sqrt{2}} \left[ \begin{array}{ccl}
0 & r(x) |x|^{\frac{1}{2}} \\ 0 & |x|^{\frac{1}{2}} \end{array} \right] , $$
we have $z z^*  = \hat{x}$. Thus $\hat{x}$ is positive.

\begin{lemma} \label{thin}  If $Z$ is a WTRO and $y \in \bar{L}(Z)$
then $y= \hat{x}$ for an $x \in {\rm Ball}(Z)$ iff $0 \leq y \leq \hat{u}$
for a tripotent $u \in Z$.  If these hold then $x \in {\mathfrak c}_u$.
\end{lemma}

\begin{proof}  If $y= \hat{x}$ then
matrix multiplication shows
that $0 \leq y \leq \widehat{r(x)}$.
For the converse, if $0 \leq y \leq \hat{u}$ then
multiplying $y$ with $\hat{u}$, we obtain that $(1/2)uu^{\ast}y_{11}+(1/2)uy_{12}^{\ast}=y_{11}$
and $(1/2)u^{\ast}y_{11}+(1/2)u^{\ast}uy_{12}^{\ast}=
y_{12}^{\ast}$. Since $y_{11} \leq (1/2)uu^{\ast}$,
it follows  that $uy_{12}^{\ast}=y_{11}=y_{12}u^{\ast}$
and $u^{\ast}y_{11}=u^{\ast}uy^{\ast}_{12}=y^{\ast}_{12}$.  Thus
it is easy to see that
$y_{11}^{2}=y_{12}y_{12}^{\ast}$. A similar argument shows
that $y_{22}^{2}=y_{12}^{\ast}y_{12}$.
The equalities above also show that
$y_{12} \in Z_{2}(u)$, and since $y_{12} u^* = y_{11} \geq 0$,
it follows that $y_{12}$ lies in ${\mathfrak c}_u$.
\end{proof}

\begin{corollary}  If $\{ u_i \} $ is an increasing net
of tripotents in a WTRO $Z$  then $u_i \to u$ weak*
iff $\widehat{u_i} \to \hat{u}$ weak*.
\end{corollary}

\begin{proof}   If $u_i \to u$ weak*,
then $\{  \widehat{u_i} \} $ is an increasing net
of projections dominated by $\hat{u}$.
Its weak* limit, by Lemma \ref{thin}, equals
$\hat{v}$ for some $v$.  Looking at convergence in
the 1-2 corner, we see $u = v$.  So $\widehat{u_i} \to \hat{u}$ weak*.
The other direction is easier.
\end{proof}

\begin{corollary} Let $\{ u_{\lambda} \}$ be a family of tripotents that
pairwise satisfy any one of the conditions in  Lemma {\rm \ref{sup}}.
Then $\widehat{\bigvee_\lambda \;  u_{\lambda}} = \bigvee_\lambda \;
\widehat{u_{\lambda}}$.   If also the $u_{\lambda}$ are all open,
then so is $\bigvee_\lambda \;  u_{\lambda}$.    \end{corollary}

\begin{proof}
 If $\hat{v} \perp \breve{u}, \hat{w} \perp \breve{u}$, and
$\hat{v} \perp \breve{w}$, then by Lemma \ref{sup} we have
$\widehat{v \vee w}=r(\hat{v}+\hat{w})$.
Since $p(\hat{v}+\hat{w}) \perp \breve{u}$
for any odd polynomial $p$, we have $\breve{u} \perp
r(\hat{v}+\hat{w}) = \widehat{v \vee w}$.
By Lemma \ref{sup}, $u \vee (v \vee w)$ exists. By induction,
 $\bigvee_{\lambda \in F} \; u_{\lambda}$
exists for any finite set $F$,
and $\widehat{\bigvee_{F} \, u_{\lambda}} = \bigvee_{F} \, \widehat{u_{\lambda}}$.
We leave the rest as an exercise, using
the last Corollary.
\end{proof}

The material in the rest of this section is used in
\cite{Betal}.

If $Z$ is a TRO
and $x \in {\rm Ball}(Z)$, then Edwards and R\"uttimann define
$u(x)$ to be the weak* limit in $Z''$, or equivalently in the
W*-algebra $Z''_2(r(x))$, of $x^{2n+1}$, where
$x^{2n+1} = x x^* x \cdots x^* x$, a product of $2n+1$ terms
(see \cite[Lemma 3.4]{ER1}).

\begin{definition} \label{snn}    If $Z$ is a TRO
then a tripotent $v$ in $Z''$ is
{\it compact} if it is the weak* limit of a decreasing net of
tripotents $u(x_{\lambda})$, where each $x_{\lambda} \in {\rm Ball}(Z)$.
\end{definition}

{\bf Remarks.}  1) \  We do not need this here, but it an easy exercise to show that
$u(x)$ is the largest tripotent $v$ such that $v = v x^* v$ (see
also \cite[Lemma 3.4]{ER1}).

2)  \ Clearly $u(x)$ is compact for any $x \in {\rm Ball}(Z)$.

3)  \ If $Z$ is a C*-algebra $A$, and $x \in A_+$, then $u(x)$ is
a projection.  Note that in this  case, if $0 \leq x \leq y \leq 1$
then $u(x) \leq x \leq y$, so that $u(x) \leq y^n$ for any $n \in \Ndb$.
Thus $u(x) \leq u(y)$.

4) \ It is essentially implicit in the main result from \cite{ER1} that
if $Z$ is a WTRO, and $x, y \in {\rm Ball}(Z)$, then
$u(x) \wedge u(y) = u(\frac{x+y}{2})$.  This is used in \cite{Betal}.

\begin{lemma} \label{thin2}  If $Z$ is a TRO and  $x \in {\rm Ball}(Z)$
then $u(\hat{x}) = \widehat{u(x)}$ and $r(\hat{x}) = \widehat{r(x)}$.
\end{lemma}

\begin{proof}  A simple computation shows that
$$\hat{x}^{2n+1} = \frac{1}{2} \left[ \begin{array}{ccl}
 |x^*|^{2n+1} & x^{2n+1} \\ (x^*)^{2n+1} & |x|^{2n+1} \end{array} \right] ,
$$
where $x^{2n+1}$ is as above Definition \ref{snn}.
The weak* limit of $\hat{x}^{2n+1}$ is a projection $q$ say, whose
$1$-$2$ entry is $\frac{1}{2} u(x)$ by the last displayed formula.
As we said earlier,  $\hat{x} \leq \widehat{r(x)}$.
This, together with Lemma \ref{thin}, shows that $q = \widehat{u(x)}$.
 That is,
$u(\hat{x}) = \widehat{u(x)}$.
If $p$ is an `odd polynomial', then
$$p(\hat{x})  = \frac{1}{2} \left[ \begin{array}{ccl}
 p(|x^*|)  & p(x) \\
p(x^*)  & p(|x|) \end{array} \right] .
$$
It follows by a norm approximation that
the same relation holds with $p$ replaced by the
function $t^{\frac{1}{2n+1}}$.
All of these quantities are bounded above by $\widehat{r(x)}$.
In the weak* limit, and using Lemma \ref{thin},
it follows that $r(\hat{x}) = \widehat{r(x)}$.
\end{proof}

Let $v$ be a tripotent in $Z''$, for a TRO $Z$.
Following \cite{AkeP},  we say that $v$
{\it belongs locally} to $Z$
if $v^{\ast}v$ is a closed projection in $(Z^* Z)''$ and
$v=xv^{\ast}v$ for an element $x \in {\rm Ball}(Z)$.
The following known result \cite{AkeP,ER4}, which we give a
quick alternative proof of for the readers convenience, shows
that this is equivalent to $v$ being compact:

\begin{proposition}  \label{iseq}   {\rm (Akemann-Pedersen, Edwards
and R\"uttimann) \ }   Let
$u$ be a tripotent in $Z''$ for a TRO $Z$.  The following are
equivalent:
\begin{itemize}
\item [(i)]  $u$ is compact.
\item [(ii)]  $\hat{u}$ is a compact projection in $L(Z)''$ (that is,
there exists a decreasing net in $L(Z)$ converging
weak* to $\hat{u}$).
\item [(iii)]  $u$ belongs locally to $Z$.
\end{itemize}
Also, a  weak* limit of a decreasing net
of compact tripotents is compact.
\end{proposition}

\begin{proof}  (i) $\Rightarrow$ (ii) \ There is a family of norm one elements $x_{\lambda}$ such
that $u(x_{\lambda})$ is a decreasing net of tripotents converging weak* to $u$.
Then $\widehat{u(x_{\lambda})}$ is a decreasing net of projections
converging weak* to a projection $p \geq \hat{u}$, say.
By Lemma \ref{thin}, we have $p = \hat{x}$
for some $x \in Z''$, and looking at the $1$-$2$ entry we see that
$x = u$.   So $\widehat{u(x_{\lambda})} \to \hat{u}$ weak*.  By Lemma \ref{thin2},
we have $\widehat{u(x_{\lambda})} = u(a)$ for some $a \in L(Z)$, and is thus
a compact projection.  Thus $\hat{u}$ is closed, being a decreasing limit
of closed projections.  Since it is bounded above
by an element in $L(Z)$, it is compact.

(ii) $\Leftrightarrow$ (iii) \   This is Proposition 4.9 of \cite{AkeP}.

(iii) $\Rightarrow$ (i) \  If
$u=xu^{\ast}u$ for an element $x \in {\rm Ball}(Z)$,
 then
 $u= uu^{\ast}x = uu^{\ast}xu^{\ast}u$ (see \cite[Lemma 4.8]{AkeP}), which
implies that
 $$x=  uu^{\ast}x u^{\ast}u + (1-uu^{\ast})x(1-u^{\ast}u) =
u + (1-uu^{\ast})x(1-u^{\ast}u) . $$
It follows that  $x^{2n+1}= u +(1-uu^{\ast})x^{2n+1}(1-u^{\ast}u)$, where
$x^{2n+1}$ as usual means $x x^* x \cdots x^* x$, a product of $2n+1$ terms.
This implies that $r(x) \geq u$.  Since $r(x)^{\ast}r(x) \geq u^{\ast} u$,
by the Urysohn lemma for C*-algebras
there is a decreasing net $(y_{\lambda})$  in $Z^{\ast}Z$
converging to $u^{\ast}u$
with $y_{\lambda} \leq r(x)^{\ast}r(x)$.
 Now  $r(x)y_{\lambda}$ lies in $Z_{2}''(r(x))$,
and so  $(r(x)y_{\lambda})$ is decreasing in $Z''_{2}(r(x))$.
Hence $x \, (r(x)y_{\lambda})^{\ast} \, x$
is a decreasing net in $Z(r(x))_{+} \cap {\rm Ball}(Z)$
converging weak* to $x u^{\ast} u r(x)^{\ast}x= uu^{\ast}x= u$.
Thus the projections $u(x(r(x)y_{\lambda})^{\ast}x)$
are a decreasing (by
Remark 3 after Definition \ref{snn}) net converging weak* to $u$.

Finally, given a decreasing net
of compact tripotents with limit $u$, a slight modification
of the first paragraph
of the proof shows that $\hat{u}$ is compact.   \end{proof}

We now give a Urysohn lemma for TRO's, based on
Akemann's Urysohn lemma for C*-algebras \cite{Ake,Ake2,Ake3,AkeP}.
See also \cite{Per} for a related result, with a different proof
strategy (which relies on results of the second author \cite{N}).

\begin{theorem}  Suppose that $Z$ is a TRO and
that $v$ and $u$ are tripotents in $Z''$
such that $v$ is compact, $u$ is open, and $v \leq u$.
Then there exists an element $x \in Z$ such that $v \leq x \leq u$ in the W*-algebra $Z''_{2}(u)$.
\end{theorem}

\begin{proof}
By Proposition \ref{prop3} we have  $\hat{v} \leq \hat{u}$.  By
 the Urysohn lemma for C*-algebras,
there is an element $y \in L(Z)_{+}$ such that $\hat{v} \leq y \leq \hat{u}$.
By Lemma \ref{thin}, we have  $y = \hat{x}$
for some $x \in {\mathfrak c}_u \subset Z''$.
By Lemma \ref{thin2} it follows that $u(y) = \widehat{u(x)}$
and $r(y) = \widehat{r(x)}$.
 By
Proposition \ref{prop3} we have  $v \leq u(x)\leq r(x) \leq u$.
Hence $v \leq x \leq u$.
\end{proof}

{\bf Remarks.}  \ 1) \ It is easy to see from this Urysohn lemma, that if $v \leq u$ are
tripotents in $Z''$ with $u$ open, then $v$ is compact in $Z''$ iff $v$ is compact
as a projection in $Z(u)'' = Z''_2(u)$.

2) \ There are also `regularity' properties for
open and compact tripotents, analoguous to the
case of projections in a C*-algebra (see e.g.\
Akemann's regularity property described in
\cite[Section 2]{H}).  The following result corresponds to the `normality' separation
property one has in locally compact topological spaces:

\medskip

The following variant of Urysohn's lemma solves an open problem from \cite{Per},
in the special case of TROs (see \cite[Problem 2.13]{Per}).

\begin{theorem}  Suppose that $Z$ is a TRO and
that $v$ and $w$ are compact tripotents in $Z''$
with $v^* w = v w^* = 0$.
Then there exist elements $x,y \in {\rm Ball}(Z)$ such that
 $r(x)^* r(y) = r(x)r(y)^* = 0$, and $v \leq x$ and $w \leq y$
in the C*-algebras $Z''_2(r(x))$ and $Z''_2(r(y))$ respectively.
\end{theorem}

\begin{proof}
Clearly $\hat{v} \perp \hat{w}$ are compact, and so
$\widehat{v+w} = \hat{v} + \hat{w}$ is a
closed projection.  Since both $\hat{v}$ and $\hat{w}$ are
dominated by an element
in $L(Z)_+$, so is $\hat{v} + \hat{w}$.  Hence $\widehat{v+w}$ is a compact
projection \cite{AkeP},
and so $v+w$ is a compact tripotent.  Thus there exists an
open tripotent $u \geq v + w$.  Working inside $Z''_2(u)$, we have
that $v, w$ are compact  mutually orthogonal projections (see Remark 1) above),
and by \cite[Proposition 2.6]{AAP} there exist mutually orthogonal
open projections $p,q$ in $Z''_2(u)$ with $v \leq p, w \leq q$.   By the
noncommutative Urysohn lemma, there exist elements $x,y \in {\rm Ball}(Z)
\cap Z''_2(u)$
such that $v \leq x \leq r(x) \leq p$ and $w \leq y \leq r(y) \leq q$,
all inequalities in the C*-algebra $Z''_2(u)$.   We leave the rest
as an exercise.
\end{proof}

One might ask if there is a Urysohn lemma for the case that one of the
tripotents is merely `closed', as in
\cite[Proposition 2.6]{AAP}.  It is clear by the methods above
that this equivalent to asking if every `closed' tripotent is dominated
by an open one, and we are not sure if the latter holds in all TROs.

 We end this section with a couple of results which
are interesting in their own right, and which we will need later.

\begin{lemma} \label{matr}   Let $A$ be a
C*-algebra and $x = [x_{ij}] \in M_n(A)_+$.
Then $$(\wedge_i \; u(x_{ii})) \otimes I_n
\leq u(x) \leq x \leq r(x) \leq (\vee_i \; r(x_{ii})) \otimes I_n .$$
\end{lemma}

\begin{proof}  By an obvious induction argument it suffices
to prove the case that $n = 2$.
Let $p = \vee_i \; r(x_{ii})$.
Clearly $x_{ii} p = x_{ii}$
for each $i$.   We claim that $x_{ij} p = x_{ij}$ for
 each $i,j$.
Note that $p^\perp r(x_{22}) p^\perp = 0$.
Thus $0 \leq (1 \oplus p^\perp) x (1 \oplus p^\perp)$,
which forces, by elementary operator theory,
that $x_{12} p^\perp = 0$.   A similar argument
shows that $x_{21} p^\perp = 0$.
Thus $x (p \otimes I_n) = x$,
so that $r(x) \leq p \otimes I_n$ as desired.

By elementary operator theory, if $q$ is a projection with
$q x q = q$ then $q x (1-q) = (1-q) x q = 0$.
Using the principle in the first paragraph of the
proof again, and the fact that
$u(x_{ii}) = x_{ii} u(x_{ii}) = u(x_{ii}) x_{ii}$
we find that $x_{ij} u(x_{ii}) = 0$ if $i \neq j$.
This gives $x ((\wedge_i \; u(x_{ii})) \otimes I_n)
= (\wedge_i \; u(x_{ii})) \otimes I_n$, so that
$(\wedge_i \; u(x_{ii})) \otimes I_n  \leq u(x)$.
\end{proof}

In the following result, for a $*$-TRO $Z$ we write $\tilde{{\mathcal L}}$
for the C*-subalgebra of $L(Z)$ whose two main
diagonal entries are equal, and whose off diagonal entries are equal (see e.g.\ \cite[Section 2]{BW}).
We call this the {\em restricted linking algebra}.

\begin{lemma}\label{proj} Let $Z$ be a $*$-TRO, and
let $p$ be a self-adjoint projection in
the center of $\tilde{{\mathcal L}}''$.
Then $p$ is of the form $\frac{1}{2}(I(u)+I(u)^{2}) \oplus^{\infty} (q
\otimes e_{11}+ q \otimes e_{22})$, where $u$ is a central selfadjoint
tripotent in $Z^{\prime\prime}$, and $q$ is a central
projection in $(Z^2)''$ such that $qu=u q=0$.  If $p$ is open in
$\tilde{{\mathcal L}}''$ then $q$ is open in $(Z^2)''$.  \end{lemma}

\begin{proof}
Let $\Theta$ as usual be the canonical $*$-automorphism on $L(Z)$,
namely changing the sign of the `off-diagonal' corners.
Then $\Theta(p)$ is central too.
Let $u = p - \Theta(p)$.  A computation shows that $u$ is a tripotent, and $\frac{1}{2}(u + u^2) =
p - p \Theta(p)$.  Since  $p$ is the orthogonal sum
of $p \Theta(p)$ and $\frac{1}{2}(u + u^2)$, it is
easy to see the first assertion.  If $p$ is open
then so is $\Theta(p)$ by an obvious argument using the
 canonical $*$-automorphism on $L(Z)$.  Since a product of
central open projections is open,
$p \Theta(p) = q \otimes e_{11}+ q \otimes e_{22}$ is open in
$\tilde{{\mathcal L}}''$,
and now it is easy to see that $q$ is open in $(Z^2)''$.
\end{proof}

{\bf Remark.} Simple examples show that in the last lemma
one cannot hope that $u$ is an open tripotent necessarily, if $p$ is open,
even if $Z$ is the commutative C*-algebra $C([0,1])$.
See, however, Lemma \ref{sees} for something along this line.

\section{Maximal cones on TROs}

We will need to develop TRO generalizations of facts from \cite[Section 5]{BW}.
The reader may wish to follow along with that paper.

If $Z$ is a TRO, and if $u$ is  a tripotent in $Z$,
then $u \otimes I_n$ is  a tripotent in $M_n(Z)$,
and hence there is an associated C*-algebra
$M_n(Z)_2(u \otimes I_n)$, which equips
$M_n(Z)$ with a cone ${\mathfrak c}_n$.  Of course
${\mathfrak c}_1 = {\mathfrak c}_u$.
Indeed, $${\mathfrak c}_n \, = \, {\mathfrak c}_{u \otimes I_n}
 \, = \,
\{ [x_{ij}] \in M_n(Z_2(u)) : [u^* x_{ij}] \geq 0 \}
 \, = \, M_n(Z_2(u))_+  .$$
Similarly, if $u$ is  a tripotent in $Z''$
then we have a canonical natural cone on $M_n(Z)$:
 $${\mathfrak d}_n =
{\mathfrak d}_{u \otimes I_n}
 \, = \, \{ [x_{ij}] \in M_n(Z(u)) : [u^* x_{ij}] \geq 0 \}
 \, = \, M_n(Z(u))_+ .$$
Sometimes we will write ${\mathfrak d}_u$ for the entire sequence
$({\mathfrak d}_{u \otimes I_n})$, and similarly for ${\mathfrak c}_u$.

We will need a fact about quotients of TROs.  First recall
that if $Z$ is a TRO, and if $J$ is a ternary ideal in $Z$,
then $Z/J$ may again be viewed as a TRO (see e.g.\ \cite[Section 8.3]{BLM}).

\begin{lemma} \label{qrs}  If $J$ is a ternary ideal in
a naturally ordered TRO $Z$, then the TRO $Z/J$ possesses a
natural cone for which the canonical
quotient ternary morphism $Z \to Z/J$ is
completely positive.
\end{lemma}

\begin{proof}
If $Z$ is a TRO in a C*-algebra $A$, we consider $Z''$ as a TRO in the
W*-algebra $A''$.
Now $J^{\perp \perp}$ is a weak* closed
ternary ideal in $Z''$,
and hence equals $Z''q$ for a central projection $q$ in $(Z^* Z)''$
as is well known (for example, it is
a special case of \cite[Theorem 7.4 (vi)]{BZ}).  If $p = 1-q$ then
$(Z/J)'' \cong Z''/J^{\perp \perp} \cong Z'' p$.
We may thus identify $Z/J$ as a TRO inside the WTRO $Z'' p$.  This endows
$Z/J$ with  natural matrix cones.  Let
$q_J : Z \to Z/J$ be the quotient ternary
morphism.
If $z \in Z_+$ then $z \geq 0$ in $A''$, and so
$z = (z^* z)^{\frac{1}{2}}$.
Thus $z p = (z^* z)^{\frac{1}{2}} p
= p (z^* z)^{\frac{1}{2}} p \geq 0$,
and so $q_J(z)$ is in the cone just defined
in $Z/J$.   A similar argument applies to matrices,
so that $q_J$ is completely positive.
\end{proof}

\begin{lemma} \label{orna}  Let $Z$ be a TRO with
matrix cones ${\mathfrak C} = ({\mathfrak C}_n)$ for which there
exists a completely positive complete isometry
from $Z$ into a C*-algebra.    Then the given
cones ${\mathfrak C}$ in $Z$ are contained in a
natural cone for $Z$.
\end{lemma}

\begin{proof}
Just as in \cite[Lemma 5.3]{BW}.
\end{proof}

\begin{definition} \label{defmax}
We say that an operator space ordering $({\mathfrak c}_n)$
on an operator space $X$
is {\em maximal}, or that $X$ is {\em maximally ordered},
if $({\mathfrak c}_n)$ is maximal amongst the operator space
orderings on $X$.  This is equivalent to saying that
every completely positive complete isometry $X \to B$ into
a C*-algebra is a complete order embedding.
\end{definition}

It follows from Lemma \ref{orna} that the maximal
(operator space) orderings on a TRO, are precisely
the maximal natural orderings.

\begin{theorem}\label{okg}
Suppose $Z$ is a TRO with
an operator space ordering.  Then $Z$ has a maximal
(operator space) ordering majorizing the given one,
and this cone is natural.
\end{theorem}

\begin{proof}
Just as in \cite[Theorem 5.4]{BW}, but including
an appeal to  Proposition \ref{incr}.
\end{proof}

As mentioned after Corollary \ref{tprop3}, natural dual cones
in a WTRO $W$ correspond bijectively to tripotents
in $W$.   This gives a very satisfactory characterization
of the maximal natural dual cones.
Maximal natural dual cones in a WTRO $W$ correspond to maximal tripotents, which are
exactly the extreme points of Ball$(W)$.  Indeed the
extreme points of Ball$(W)$ are well known to be the tripotents
such that $(1- u u^*) W (1 - u^* u) = (0)$.  Any such tripotent is
maximal, since if $v \geq u$ then $$0 = (1- u u^*) v (1 - u^* u)
= (v - u) (1 - u^* u) = v - u .$$
Conversely,
if the WTRO $(1- u u^*) W (1 - u^* u)$ is not $(0)$
then it has a nonzero tripotent $w \perp u$, and $w + u \geq u$.  Thus $u$ is
not maximal.

In a TRO $Z$, maximal natural  cones correspond to maximal open
tripotents.   These exist by Zorn's lemma, since any increasing
chain of open tripotents is bounded above by an open tripotent
(Proposition \ref{incr}).    We consider maximal open tripotents in
Theorem \ref{gentr} below, also settling an issue raised in
\cite[Section 5]{BW}, which we now describe.

Let $A$ be a C*-algebra and let $p$ and $q$ be an open and a closed central
projection in $A^{\prime\prime}$.  We say that $q$ {\it is contained in the
boundary} of $p$, if $p \perp q$ and if whenever $r$ is an open
central\footnote{Dealing with {\em central} projections here yields a
simpler characterization without introducing any additional complications.}
projection in $A^{\prime\prime}$ which is
perpendicular to $p$, then $r$ is perpendicular to $q$.
We shall not use this, but if
$q + p$ is closed, which will be the case for us below,
then it is easy to see that $q$  is contained in the
boundary of $p$ iff $q + p$ is the smallest closed
central projection dominating $p$.
 In \cite[Proposition 5.11]{BW} it was shown that  if $Z$ is a
$*$-TRO and $u$ is a selfadjoint central open tripotent in
$Z''$, then $u$ is maximal amongst the
selfadjoint central open tripotents in $Z''$ if
$1-u^2$ is contained in the
boundary\footnote{We shall not use this, but it is easy to see
that $1-u^2$ is contained in the
boundary of $\hat{u} = \frac{1}{2}(u + u^2)$ iff
$\frac{1}{2}(u^2-u)^{\perp}$ is the smallest closed
central projection dominating $\frac{1}{2}(u + u^2)$.}
of both $\frac{1}{2}(u + u^2)$ and $\frac{1}{2}(-u + u^2)$.
Here the C*-algebra $A$ is $Z + Z^2$.
It was also noted there that the converse of this  is true in the
`commutative
case'; however we remark that an inspection of the
proof of this converse (see \cite[Corollary 6.8]{BW})
shows that we were also assuming there that
$Z \cap Z^2 = (0)$.  This is not a serious
restriction, since any $*$-TRO is ternary $*$-isomorphic to
one satisfying this property, and below we shall always assume that
$Z \cap Z^2 = (0)$ when we use the phrase
`$1-u^2$ is contained in the boundary of $\frac{1}{2}(u + u^2)$'.
 It was suggested in \cite{BW} that such a
`contained in the boundary' condition might characterize
maximal selfadjoint central open tripotents for any
$*$-TRO.  By results in
that paper, such a characterization would immediately
give a characterization of the maximal ordered operator space cones,
and thus also the maximal  cones which are natural in the sense of
\cite{BW}, on
any $*$-TRO.

To motivate the value of having such a characterization, it is very
instructive to look at a commutative example studied
in \cite[Section 6]{BW}.  Let $S^2$ be the
unit sphere, and $Z$ the $*$-TRO
$\{ f \in C(S^2) : f(-x) = - f(x) \}$.   In this case
open selfadjoint tripotents $u$ in $Z''$ correspond precisely
to open subsets $U$  of the sphere (called blue), which do not intersect
$-U$ (called red).   Suppose that $S^2 \setminus (U \cup (-U))$
is colored black.  The `contained in the
boundary' characterization discussed in the
last paragraph says precisely\footnote{This is misstated
in the fourth last line of p.\ 709 of \cite{BW}, but
the typographical error should have been clear in
the context.} that $u$
(and hence the associated ordering of $Z$)
is maximal iff the black region is the boundary
of the red region (and hence also of the blue region).
Thus, for example,  a sphere whose top hemisphere is red
and whose bottom hemisphere is blue, with a black equator line,
is maximal; whereas if one were to thicken the
equator to a black band one loses maximality.
From the geometry of such examples, it seems clear that one could not
improve on this characterization.

In the `noncommutative case' it is unfortunately
not true that if $u$ is maximal amongst the
selfadjoint central open tripotents in $Z''$
then $1-u^2$ is contained in the
boundary of $\frac{1}{2}(u + u^2)$.
For example, take $Z$ to be the subspace of $M_2$ with
`main diagonal' entries zero.  In this case, $u = 0$ is a
maximal  selfadjoint central open tripotent in $Z''$,
but $r = I_2$ satisfies $r \frac{u + u^2}{2} = 0$
but $r(1-u^2) = r \neq 0$.    It does not help if we
replace $A = Z + Z^2$ by the `restricted linking C*-algebra'
$\tilde{{\mathcal L}}$ mentioned at the end of Section 4, indeed this situation is equivalent
since $Z + Z^2$ is $*$-isomorphic to $\tilde{{\mathcal L}}$ if
$Z \cap Z^2 = (0)$.

This problem can be remedied in
 several ways.  For example, we can put
a restriction on the projections $r$ considered in the
definition of `contained in the
boundary'.  If $A=Z+Z^{2}$ and $Z \cap Z^{2}=0$, we say that a
projection $r$ in $A^{\prime\prime}$ is
{\em antisymmetric}
if $\Theta^{\prime\prime}(r) \perp r$.  Here $\Theta:A \rightarrow A$ is
the period 2 $*$-automorphism
$\Theta(z+a)=a-z$
for $z\in Z$ and $a \in
Z^{2}$.  If $r$ is a
projection (resp. central projection) in $A^{\prime\prime}$, then
it is easy to see that $r$ is antisymmetric
iff $r = \frac{v + v^2}{2}$ for a selfadjoint tripotent (resp. central
selfadjoint tripotent) $v \in Z''$.
If $r$ is an open central projection in $A^{\prime\prime}$, then
it is easy to check using Lemma \ref{proj} and \cite[Proposition 4.18]{BW},
that $r$ is antisymmetric
iff $r$ dominates no nontrivial open central projection in $(Z^2)''$.
We say that $q$ is {\em antisymmetrically contained in the
boundary} of $p$,
if $p \perp q$, and whenever $r$ is an open
antisymmetric central projection in $A^{\prime\prime}$ which is
perpendicular to $p$, then $r \perp q$.

For a general TRO $Z$, we use the definition of
antisymmetric projections from the introduction.  Let $A = L(Z)$,
and let $p$ and $q$ be respectively open and  closed projections in $A^{\prime\prime}$.
We say that $q$ is {\it antisymmetrically contained in the boundary} of $p$, if
$p \perp q$, and whenever $r$ is an open antisymmetric
projection in $A^{\prime\prime}$ which is perpendicular to $p$, then $r \perp q$.

The following is a characterization
of maximal open tripotents in a TRO:

\begin{theorem} \label{gentr} Let $Z$ be a TRO.
Suppose $u$ is an open tripotent in  $Z^{\prime\prime}$. Then $u$ is maximal amongst
the open tripotents if and only if $1-I(u)^{2}$ is antisymmetrically
contained in the boundary of $\hat{u}$.  \end{theorem}

\begin{proof}
Since $u$ is maximal iff $-u$ is maximal, we may replace $\hat{u}$ by $\breve{u}$.
Suppose that  $1-I(u)^{2}$ is antisymmetrically
contained in the boundary of $\breve{u}$. If $v \geq u$ and if $v$ is
open, then $v$ commutes with $u$, and $r = \hat{v}$ is (by Theorem \ref{chopt})
an open projection in $L(Z)''$.  It is easy to check that
$r \breve{u} = 0$.  Thus $r (1-I(u)^{2}) = 0$ which gives
$v = v u^* u = u$.  Thus $u$ is maximal.

Conversely, suppose that
$u$ is a maximal open tripotent in $Z^{\prime\prime}$.
Suppose that  $r=\hat{v}$ is an open antisymmetric projection with
$r \breve{u} = 0$.  From the ensuing commutator relations
coming from this last equality, we see that
$u \vee v$ exists by Lemma \ref{sup}, and thus it is open
by Proposition \ref{prop4}.  Since $u$ is maximal,
$u = u \vee v$, and so $v \leq u$.  Hence $v = vu^{\ast}u$.
Inspection now reveals that $\hat{v} (1-I(u)^{2})=0$.
\end{proof}

The following result, whose proof we omit since
it is essentially the same as
the proof of Theorem \ref{gentr},
characterizes maximal cones in a $*$-TRO.

\begin{proposition} \label{newl} Let $Z$ be a $*$-TRO, and let $u$ be a selfadjoint
central open
tripotent in  $Z^{\prime\prime}$.  Then $u$ is maximal amongst the
selfadjoint central open tripotents if and only if $1-u^{2}$ is
antisymmetrically contained in the boundary of $\frac{u + u^2}{2}$.  \end{proposition}

{\bf Remark.}  It follows easily from what we have done, that
the conjecture from \cite{BW} that we have been discussing,
is true for the class of TROs $Z$ which have the
following property: whenever $p$ and $q$ are respectively open and closed
central projections in $(Z + Z^2)''$, which are not orthogonal to each other,
then $p$ dominates an antisymmetric open central projection $r$
which is not orthogonal to $q$.    This condition is always
satisfied in the `commutative case' of \cite[Section 6]{BW}.
We remark too that the fact that the conjecture is true
in this `commutative case',
also follows from the following result (since
in the commutative case, in the notation below, $r$ is necessarily
antisymmetric, for if it were not
then the existence of $v$ below contradicts the maximality of $u$).

\begin{lemma} \label{sees}
Let $Z$ be a $*$-TRO with  $Z \cap Z^{2} = (0)$.
Suppose that $r$ is an open central projection in
$(Z + Z^2)''$ which
is perpendicular to $\hat{u} = \frac{1}{2} (u + u^{2})$
for a maximal central selfadjoint open tripotent $u \in Z''$.
Then either $r$
is antisymmetric (and thus orthogonal to $1 - u^{2})$,
or there exists a nonzero selfadjoint open tripotent $v \in Z''$
 which is perpendicular to $u$ and $r \geq \hat{v}$.
\end{lemma}

\begin{proof}
By  Lemma \ref{proj}, $r = q + \frac{1}{2}(u+u^{2})$ where $q$ is an open
central  projection in $Z^{2}$, $u$ is a central selfadjoint
tripotent, and $q \perp u$. If $q = 0$ then we are done:
$r$ is antisymmetric and is
orthogonal to $1 - u^{2}$ by Proposition \ref{newl}.
If not, suppose that $x_{\lambda}$
is a net in $Z^{2}$ converging up to $q$. Clearly
there exists an element $y \in Z$, which we
can take to be selfadjoint, and a $\lambda$ with
$z=x_{\lambda} y x_{\lambda} \ne 0$ (for otherwise,
taking a strong limit in $(Z + Z^2)''$, we have
$q y q = q y = 0$ for all $y \in Z$, so that $q = 0$).
Of course $r(z)$ is open in $Z''$.
We have  $z = zq \perp u$, and so $r(z) \perp u$ and $r(z) q = r(z)$.
Clearly
$\widehat{r(z)} \leq q \leq r$.
 \end{proof}

Let $Z$ be a $*$-TRO.  For any set $S \subset Z''$, we denote by $S_\vdash$
the set $\{x \in Z : yx=xy=0 \mbox{ }\forall y \in S\}$.
If ${\mathfrak c}$ is a natural dual cone in $Z''$ in
the sense of \cite{BW},
 then ${\mathfrak c} = {\mathfrak c}_u$ for an open central
selfadjoint tripotent $u \in Z''$.   In this case,
${\mathfrak c}_\vdash = u_\vdash$, and this is a ternary $*$-ideal
in $Z$.    It follows from \cite[Lemma 3.4]{BW} that
$u_\vdash = \{ z \in Z : r z = z \}$, and $(u_\vdash)^{\perp \perp}
= r Z''$, for an open central projection in $(Z^2)''$ such that
$r z = z p$ for all $z \in Z$.   We claim that $r u = 0$.
Indeed $r u \in r Z'' = (u_\vdash)^{\perp \perp}$,
and if $x_t \in u_\vdash$ with $x_t \to r u$ weak*, then
$0 = u^2 x_t \to r u$, so that $r u = 0$.

\begin{corollary}  \label{insmat}  If $Z$ is a $*$-TRO,
and $u$ is an open selfadjoint tripotent in $Z''$ then
$1-u^2$ is contained in the boundary of $\frac{1}{2}(u+u^{2})$
iff $u$ is maximal amongst the
open selfadjoint central tripotents and
$u_\vdash = (0)$.
\end{corollary}

\begin{proof}   ($\Rightarrow$) \ Under this hypothesis,
we have that $u$ is maximal as before.  If $u_\vdash \neq 0$,
then if $r$ is as above Corollary \ref{insmat} then
$r \neq 0$.  On the other hand,
$r \frac{1}{2}(u+u^{2}) = 0$, and so $r(1 - u^2) = r = 0$,
a contradiction.

($\Leftarrow$) \
If $v$ is open and perpendicular to $u$, and if
$x_t \to v$ weak* with $x_t \in Z(v)$, then $x_t \in u_\vdash = (0)$.
Thus $v = 0$.   The result then follows from Lemma \ref{sees}.
   \end{proof}

We now isolate the class of $*$-TROs for which the
conjecture from \cite{BW} is correct.

\begin{definition} \label{defcu}
A $*$-TRO $Z$ is said to be {\it completely orderable}
if for every natural dual cone ${\mathfrak c}$ in $Z''$,
either ${\mathfrak c}_{\vdash} = (0)$ or
${\mathfrak c}_{\vdash}$ has a nontrivial natural ordering
in the sense of \cite{BW}.
\end{definition}

\begin{lemma} \label{newm}  A $*$-TRO $Z$ is completely orderable iff
for every maximal open selfadjoint central tripotent $u \in Z''$ we have
$u_{\vdash} =(0)$.
 \end{lemma}

\begin{proof}  We may
assume without loss of generality
that $Z \cap Z^{2} = (0)$, and we denote $Z +Z^{2}$ by $A$.

Suppose that $u_{\vdash} =(0)$ for every maximal open selfadjoint
central tripotent $u$.  If $Z$ is not completely orderable, then
there is an open central selfadjoint tripotent $u$ such that
$u_{\vdash}$ is not orderable (and nontrivial). Write $u_{\vdash} =
\{ z \in Z : r z = z \}$ as above Corollary \ref{insmat}.   Let $v$
be a maximal central open tripotent with $v \geq u$. Since by
hypothesis $v_{\vdash} = (0)$, we must have $vr \ne 0$ (for if $vr =
0$ then if $0 \neq x \in u_{\vdash}$ then $x v = x r v = 0 = v r x =
v x$, so that $0 \neq x \in v_{\vdash} = (0)$). Let  $w=vr$, a
central selfadjoint tripotent in $Z''$ which is perpendicular to $u$
by the line above Corollary \ref{insmat}.  Note that $r$ is also an
open central projection in $A''$, and that $v$ is an open
selfadjoint central tripotent in $A''$, so that $w$ is an open
selfadjoint central tripotent in $A''$ by e.g.\ Corollary
\ref{prop1}. Suppose that $(x_{\lambda}+y_{\lambda})$ is an
increasing net in $A(w)_{+}$ which converges to $w$, with
$x_{\lambda} \in Z^{2}$ and $y_{\lambda} \in Z$.  Let $\Theta$ be
the map mentioned above Theorem \ref{gentr}. Since $-\Theta$ is a
ternary isomorphism on $A$ and $-\Theta(w)=w$, we have $-
\theta(A(w)_{+}) = A(w)_{+}$ by Proposition \ref{imop}. Hence,
$-\Theta(x_{\lambda}+y_{\lambda}) = -x_{\lambda}+y_{\lambda}$ is
also an increasing net with limit $w$ in $A(w)_{+}$. Thus
$y_{\lambda} \in A(w)_{+} \subset A''_2(w)_+$. Since $y_{\lambda}
\in Z''_2(w)$, and since $Z''_2(w)$ is a C*-subalgebra of
$A''_2(w)$, we have that $y_{\lambda} \in Z''_2(w)_+ \cap Z =
{\mathfrak d}_w$. Since  $y_{\lambda} \to w$ weak* we conclude that
 $w$ is open in $Z''$.
Since $w$ is perpendicular to $u$ we
must have $Z(w) \subset u_{\vdash}$,
and so $w$ is a tripotent in $(u_{\vdash})^{\perp \perp}$
which is open in that space.
By \cite[Corollary 5.7]{BW} we see that $u_{\vdash}$ is orderable.
This is a contradiction, and so $Z$ is completely orderable.

For the other direction, suppose that $Z$ is completely orderable.
If $u$ is a maximal open  selfadjoint central tripotent in $Z''$
 with $u_{\vdash} \ne (0)$,
then $u_{\vdash}$ has a nontrivial natural ordering.
Thus by \cite{BW}, there is
a nontrivial open  selfadjoint central tripotent $w \in (u_{\vdash})''  \cong
(u_{\vdash})^{\perp \perp}$.   This
tripotent $w$
 is also an open  selfadjoint tripotent in $Z''$.  Since
 $(u_{\vdash})^{\perp \perp} = r Z''$ as above
Corollary \ref{insmat}, and since $w r^\perp = 0$,
it is easy to see that $w$ is central in $Z''$.
Since $w \in (u_{\vdash})^{\perp \perp}$, we have
$w \perp u$.  Since $u + w \geq u$
we have arrived at a contradiction.   Thus
 $u_{\vdash}  = (0)$.
\end{proof}

\begin{theorem} \label{mtd}  If $Z$ is a completely orderable $*$-TRO
with $Z \cap Z^{2}=(0)$,
then an open selfadjoint central tripotent $u \in Z''$ is maximal
amongst the open  selfadjoint central tripotents iff $1-u^{2}$
is contained in the boundary of $\frac{1}{2}(u+u^{2})$. If $Z$ is not completely orderable,
then there exists a maximal open selfadjoint central
tripotent  $u$ such that $1-u^{2}$ is not contained in the boundary of
$\frac{1}{2}(u+u^{2})$.
\end{theorem}

\begin{proof}
If $Z$ is completely orderable,
and if $u$ is a maximal open  selfadjoint central tripotent in $Z''$,
then $u_{\vdash} = (0)$,
and so $1-u^{2}$
is contained in the boundary of $\frac{1}{2}(u+u^{2})$ by Corollary  \ref{insmat}.
The converse direction also follows from Corollary  \ref{insmat}.

Suppose that $Z$ is not completely orderable. By
 Lemma \ref{newm},
there is a maximal open selfadjoint central tripotent $u$, with
$u_{\vdash} \ne (0)$. By the lines above Corollary \ref{insmat},
there is a nonzero open central projection $r \in (Z^{2})''$
with $r \perp u$.  Thus $r$ is an open central projection in $A''$,
and $r \perp \frac{1}{2}(u+u^{2})$ but $r(1-u^{2})=r \ne 0$.
\end{proof}

{\bf Remark.}  Commutative $*$-TRO's are completely orderable,
by results in \cite[Section 6]{BW} or by a simple direct argument.
C*-algebras also satisfy this condition. In this case, for any
maximal open selfadjoint central tripotent $u$, we have that $u_{\vdash}$
is a two-sided ideal.
If $p$ is the support projection of this ideal, then
$p \perp u$.   Thus $u + p$ is an open selfadjoint central tripotent
dominating $u$.  By maximality of $u$ we have $p = 0$
and $u_{\vdash}=(0)$.

\section{Cones on operator spaces and the Shilov boundary}

In this section we study an operator space $X$ with a
given cone ${\mathfrak c}$, or with a  sequence of
matrix cones ${\mathfrak c}_n \subset M_n(X)$.  One of the
advantages of our approach is that it can be done in either
of these two settings, that is, for nonmatricial or for matricial cones.
Nonetheless we will usually focus on the matricial cone case, leaving
the nonmatricial case to the reader,
with the following lemma being an exception:

\begin{proposition}  \label{BK}  Suppose that $X$ is an operator space
with a cone ${\mathfrak c}$ which densely spans $X$, and that
$i : X \to B$ is a positive complete isometry
from $X$ into a C*-algebra.  Then the TRO $W$ generated by
$i(X)$ is a C*-subalgebra of $B$, and $W''$ is a
W*-subalgebra of $B''$.  Moreover, if $u$ is
the tripotent associated with the natural cone of $W$,
then $W = W(u)$.
\end{proposition}  \begin{proof}
If $x \in {\mathfrak c}$ then $i(x) \in W \cap B_+ \subset J(W)$
(see Lemma \ref{cruc}).
Hence $i(X) \subset J(W)$,
so that $W \subset J(W)$.  Thus $W = J(W)$ is a C*-subalgebra of $B$.
 The other assertions are obvious.
  \end{proof}

\begin{lemma} \label{bailed}
Let $({\mathfrak c}_n)$ be an operator space ordering
on an operator space $X$, and let $i : X \to B$ be a
completely positive complete isometry into
a C*-algebra.   If $W = \langle i(X) \rangle$,
set $u = \vee_{x \in {\mathfrak c}_1} \; r(i(x))$, an open
tripotent in $W''$.
Then $i_n({\mathfrak c}_n) \subset {\mathfrak d}_{u \otimes I_n}$.
\end{lemma}

\begin{proof}   We need to show that if
$[x_{ij}] \in {\mathfrak c}_n$,
and $x = [i(x_{ij})]$, then $x \in {\mathfrak d}_{u \otimes I_n}$.
  By Lemma \ref{toadd}  this
is equivalent to saying that $r(x) \leq u \otimes I_n$, which
in turn follows from Lemma \ref{matr}.
\end{proof}

As we said early in Section 4, the cone ${\mathfrak d}_{u \otimes I_n}$
in the last lemma is the natural cone in $M_n(Z)$
corresponding to the tripotent $u$ (that is, it is the
$n$'th cone in the sequence of matrix cones associated
with $u$).

Next,  we construct an ordered version of
the `noncommutative Shilov boundary' or `ternary envelope'
\cite{Ham,BLM}.   We recall its universal property,
which we use frequently.   The ternary envelope
of an operator space  $X$ is a pair $({\mathcal T}(X),j)$ consisting of
a TRO ${\mathcal T}(X)$ and a  completely isometric linear map $j : X
\to {\mathcal T}(X)$, such that
${\mathcal T}(X)$ is generated by $j(X)$ as a TRO (that is, there is no closed
subTRO containing $j(X)$), and which has the following
property: given any completely isometric
linear map $i$ from $X$ into a TRO $Z$ which is
generated by $i(X)$, there exists a (necessarily unique and surjective)
ternary morphism $\theta : Z \to {\mathcal T}(X)$ such that
$\theta \circ i = j$.
A pair $({\mathcal T}(X),j)$
with this universal property is unique up to ternary isomorphism
`fixing the copy of $X$'.    By considering simple examples
(for example, orderings on $\Cdb$!), one quickly sees that
if one wants an ordered version of this that works for operator spaces
with sensible positive cones,
the embeddings $i : X \to Z$ occurring in the universal property
above cannot be
allowed to be arbitrary completely positive complete isometries,
or even
arbitrary completely isometric complete order embeddings (unless
we have a strong extra condition,
such as $X_+$ densely spanning $X$).
We will usually need to limit the size of the cone of $Z$.

More specifically, suppose that $X$ is an operator space possessing a
cone ${\mathfrak c}$ (resp.\ sequence of matrix cones ${\mathfrak c}
= ({\mathfrak c}_n)$) such that
there is a positive (resp.\ completely positive) complete isometry
$i : X \to B$ into a C*-algebra $B$.
Then we can assign a canonical cone (resp.\ sequence of matrix cones) to
the ternary envelope $({\mathcal T}(X),j)$, namely
the intersection of all natural cones containing
$j({\mathfrak c})$ (resp.\ $(j({\mathfrak c}_n))$).
We call ${\mathcal T}(X)$ equipped with this cone structure the
{\em ordered ternary envelope} ${\mathcal T}^o(X)$.
To see that there exists at least one such cone, note that
if $i : X \to B$ is as above,
and if $W$ is the TRO generated by
$i(X)$, then by the  universal property of the
ternary envelope above, there is a
ternary morphism $\theta : W \to {\mathcal T}(X)$
with $\theta \circ i = j$.  Thus ${\mathcal T}(X)$ is
ternary isomorphic to a quotient of $W$.  By Lemma \ref{qrs},
this quotient of $W$
has a natural  cone containing the image of $i({\mathfrak c})$
in the quotient.
Hence ${\mathcal T}(X)$ has a natural cone containing  $j({\mathfrak c})$
(resp.\ containing  the sequence $(j_n({\mathfrak c}_n))$).
In particular,  $j : X \to {\mathcal T}^o(X)$ is positive (resp.\ completely positive).
It is easy to see (using Lemma \ref{bailed} in the
`matricial cone case'), that the open
tripotent corresponding to ordering which we have given
${\mathcal T}^o(X)$ is $u = \vee_{x \in {\mathfrak c}} \; r(j(x))$.

If $({\mathcal T}^o(X),j)$ is the ordered
ternary envelope of $(X,{\mathfrak c})$,
then we define the {\em order completion}
of ${\mathfrak c}$ to be the cone $\overline{{\mathfrak c}}
= j^{-1}({\mathcal T}^o(X)_+ \cap j(X))$ in $X$.
It is of interest to know
when ${\mathfrak c}$ is {\em complete}, that is,
${\mathfrak c}  = \overline{{\mathfrak c}}$,
or equivalently, that the canonical embedding of $X$ in
${\mathcal T}^o(X)$
is a (complete) order embedding.  Later in this section we
will give some sufficient conditions for this.

The following theorem is stated in the `matricial cone case'; in the
`nonmatricial case' delete the occurrences of the word `completely',
and ignore matrix cones.

\begin{theorem} \label{maintro2}  Suppose that $X$ is an operator space
with matrix cones ${\mathfrak c} = ({\mathfrak c}_n)$, and that
$i : X \to B$ is a completely positive complete isometry
from $X$ into a C*-algebra, such that if
$W$ is the TRO generated by
$i(X)$ then there is
no smaller natural cone
than $W \cap B_{+}$ on  $W$ which contains $i({\mathfrak c})$
(or, equivalently, that $W \cap B_{+}$ is the intersection of
natural cones containing $i({\mathfrak c})$).
Let $({\mathcal T}^o(X),j)$ be the ordered ternary envelope of $X$.
Then there exists
a completely positive ternary morphism $\theta : W \to {\mathcal T}^o(X)$
such that $\theta \circ i = j$.
Moreover, $\theta$ restricts to  a surjective $*$-homomorphism between
the C*-algebras associated with the natural orderings.
In particular, $\theta(W \cap B_+) = {\mathcal T}^o(X)_+$.
  \end{theorem}

\begin{proof}  Let ${\mathfrak c}$ be the cone on $X$.
By the universal property of ${\mathcal T}(X)$,
there exists such a map $\theta$; we
need to prove that $\theta$ is
completely positive.  Note that all of the `positivity' in
$Z$ `resides' in the W*-algebras $M_n(Z''_2(u))$, where $u$
is the open tripotent giving the ordering on $W$.
Similarly for
${\mathcal T}(X)$, and write $w$
 for the open tripotent giving its ordering.
 Let $\pi$ be $\theta''$ restricted to
$Z''_2(u)$, which is a weak* continuous
$*$-homomorphism (and therefore
automatically completely positive)
from $Z''_2(u)$ onto $Z''_2(\pi(u))$.
Products below are taken in those algebras.
  It suffices to show that $\pi(u) = w$.
Observe that
$\pi(i(x)^{2n-1}) = j(x)^{2n-1}$ for any $x \in {\mathfrak c}
\cap {\rm Ball}(X)$ and  $n \in \Ndb$.
Using the fact stated after
Lemma \ref{Matt},
we see that $\pi(i(x)^{\frac{1}{2n-1}}) = j(x)^{\frac{1}{2n-1}}$,
and in the weak* limit, $\pi(r(i(x))) = r(j(x))$.
Taking suprema, since weak* continuous $*$-homomorphisms
preserve suprema, we know that $\pi(u)$ is the supremum in
the W*-algebra $Z''_2(\pi(u))$ of the projections $r(j(x))$
in that algebra.  By definition of $w$, we have
$w \leq \pi(u)$.  Hence $w$ is a projection in $Z''_2(\pi(u))$,
and now it is clear that $w = \pi(u)$.  The desired surjectivity
follows from Proposition \ref{imop}.
\end{proof}

{\bf Remark.}   The ordering we have given to
${\mathcal T}(X)$ does not depend on the particular
ternary envelope chosen.  This follows immediately for example
from the universal property in the theorem.

\medskip

The ordered noncommutative Shilov boundary is particularly nice in
the case that $X$ has a densely spanning cone, for example this
boundary is a C*-algebra.  In this case, it is easy to see that we
may assume that $X$ is a selfadjoint operator space, and then the
following result is in \cite{BNW}. We include an alternative proof:

\begin{corollary}  \label{BK22}  Suppose that $X$ is an operator space
with a cone ${\mathfrak c}$ which densely spans $X$, and that
$i : X \to B$ is a positive complete isometry
from $X$ into a C*-algebra.  Then the TRO $A$ generated by
$i(X)$ equals the C*-subalgebra of $B$ generated by $i(X)$,
 and the hypothesis
on the cone of $A$ in Theorem {\rm \ref{maintro2}} holds
automatically.  Moreover, the
ordered ternary envelope of $X$ is a C*-algebra, and
the canonical ternary morphism $\theta : A \to {\mathcal T}^o(X)$
such that $\theta \circ i = j$, is a $*$-homomorphism.
\end{corollary}

\begin{proof}
The first assertions follow from
Proposition \ref{BK}.    Write $(D,j)$ for the
ordered ternary envelope of $X$, viewed as a C*-algebra.
By the universal property of the ternary envelope, there exists a
surjective ternary morphism  $\theta : A \to D$
with $\theta \circ i = j$.
Let ${\mathfrak d}$ be the intersection of the
natural cones containing $i({\mathfrak c})$.  This is
 a natural cone in $A$, and its span is an inner ideal
$J$ of $A$.  Since $J$ is a subTRO too,
$J$ contains the subTRO generated by
$i({\mathfrak c})$.   Since ${\mathfrak c}$ densely spans $X$,
$J$ contains the subTRO generated by
$i(X)$.  So $J$ = $A$, and
it follows that ${\mathfrak d} = A_+$.  Hence, and by
the `nonmatricial case' of  Theorem {\rm \ref{maintro2}},
$\theta$  is positive on $A$.   By Lemma \ref{trp},
$\theta$ is a $*$-homomorphism. \end{proof}

\begin{theorem}  \label{centl} Suppose that $X$ is an operator space
with matrix cones ${\mathfrak c} = ({\mathfrak c}_n)$,
and let  $j : X \to {\mathcal T}(X)$
be the canonical `Shilov boundary embedding'.  Then there exists
a completely positive complete isometry
from $X$ into a C*-algebra, if and only if
${\mathfrak c} \subset j^{-1}({\mathfrak d}_u)$, where
$u$ is an open tripotent in ${\mathcal T}(X)''$.

If these hold and if ${\mathfrak c}_1$ densely spans $X$, then
${\mathfrak c} = j^{-1}({\mathfrak d}_u)$ for
some open tripotent $u$ in ${\mathcal T}(X)''$ if and only if
 $(X,{\mathfrak c})$ is maximally ordered.
\end{theorem}

\begin{proof}  The ($\Leftarrow$) directions
are easy,
by looking at $j : X \to {\mathcal T}(X)$; equipping the latter space
in the proof of the first `iff'
with  the natural cone ${\mathfrak d}_u$,
and in the second `iff' with  the natural cone of ${\mathcal T}^o(X)$.

For the ($\Rightarrow$) direction of the first `iff',
let $i : X \to B$ be a completely positive complete isometry
from $X$ into a C*-algebra $B$, and let $W = \langle i(X) \rangle$,
the TRO generated by $i(X)$.
Endow $W$ with the smallest natural cone
which contains $i({\mathfrak c})$.
Let $({\mathcal T}^o(X),j)$ be the ordered ternary envelope
of $X$, and let $u$ be the
open tripotent discussed above the last theorem.
By that result, there exists
a completely positive ternary morphism $\theta : W \to {\mathcal T}^o(X)$
such that $\theta \circ i = j$.  This implies that
$j({\mathfrak c}) \subset {\mathfrak d}_u$.

Finally, if ${\mathfrak c} = j^{-1}({\mathfrak d}_u)$ for some open tripotent,
then this holds with $u$ the open tripotent corresponding to the
natural cone of ${\mathcal T}^o(X)$.   That is, ${\mathfrak c} =
\overline{{\mathfrak c}}$.  If $i(x) \geq 0$ for
a completely positive complete isometry $i$
from $X$ into a C*-algebra, then by Corollary \ref{BK22}
we have $j(x) \geq 0$, so that $x \in \overline{{\mathfrak c}} = {\mathfrak c}$.
So $i$ is an order embedding, and similarly it is a
complete order embedding.   \end{proof}

{\bf Remark.}
The `nonmatricial cone' case of the last result is valid
with the same proof.  Thus  if $X$ is an operator space
with a cone ${\mathfrak c}$, then there exists
a positive complete isometry
from $X$ into a C*-algebra, if and only if
${\mathfrak c} \subset {\mathfrak d}_u$, where
$u$ is an open tripotent in ${\mathcal T}(X)''$.

\medskip

The previous results have
nice consequences concerning unitizations, which are
explored a bit further in \cite{BNW}.
If $(X,{\mathfrak c})$ is an operator space
with a densely spanning operator space cone,
let $A = {\mathcal T}^o(X)$ be its ordered
ternary envelope, which we now know is a C*-algebra.
Let $X^1$ be the span of $X$ and the identity of the
$C^*$-algebra unitization of this $C^*$-algebra.
Then if $H$
is a Hilbert space, and
$i : X \to B(H)$ is a
completely positive complete isometry,
then it is easy to show, from the universal property of
${\mathcal T}^o(X)$, that there is a completely positive unital map from
$i(X) +  \Cdb I_H   \to   X^1$
extending the canonical map $i(X) \to X$.
One clearly has the following rigidity result:
a unital completely positive linear map $\Phi : X^1 \to B(H)$
is a complete order embedding if its restriction to
$X$ is a completely isometric complete order embedding.

\begin{corollary}  \label{fma2}
Let $X$ be an operator space with an
operator space  cone which densely spans $X$.  The following are equivalent:
\begin{itemize}
\item [(i)] The embedding of
$X$ in the unitization $X^1$  is a complete order embedding.
\item [(ii)] $X$ is maximally ordered.
\item [(iii)] The  cone ${\mathfrak c}$ is complete
(that is, the canonical embedding of $X$ in ${\mathcal T}^o(X)$ is a
complete order embedding).
\end{itemize}
\end{corollary}

\begin{proof}  We already saw in Theorem \ref{centl} that
(ii) $\iff$ (iii), and the equivalence of (i) and (iii) is obvious
from the definition of $X^1$.
\end{proof}

For the next result we recall that a {\em unital operator space}
is an operator space $X$ for which  there exists a
linear complete isometry $\varphi : X \to A$ into a
unital C*-algebra with $1_A \in \varphi(X)$.  We will
write $1$ for $\varphi^{-1}(1_A)$.  Any
unital operator space has a {\em C*-envelope} $(C^*_e(X),j)$
(see e.g.\ \cite[Section 4.3]{BLM}) which is a unital
C*-algebra together with a complete isometry
$j : X \to C^*_e(X)$ with $j(1) = 1$ such that
$j(X)$ generates $C^*_e(X)$ as a C*-algebra, and possessing
a certain universal property spelled out in the last reference.
This C*-algebra obviously has a canonical cone $C^*_e(X)_+$.

\begin{corollary}  \label{DDz}  If $X$ is a
unital operator space and if $X$ also has a
operator space cone ${\mathfrak c}$
which densely spans $X$, and which
contains $1$, then the C*-envelope $C^*_e(X)$
(in the sense of {\rm \cite[Section 4.3]{BLM}}) is
also the ordered ternary envelope of $X$, with
$C^*_e(X)_+ = {\mathcal T}^o(X)_+$.
\end{corollary}

\begin{proof}   It is known that the unital
C*-algebra $A = C^*_e(X)$ is a ternary envelope of $X$,
and hence there does exist a natural cone ${\mathfrak d}$
in $A$ containing $j({\mathfrak c})$, and therefore
also containing $1$.   The span $J$
of ${\mathfrak d}$ is a subTRO of $A$ containing $j({\mathfrak c})$,
and therefore containing $j(X)$, and $1$.  This forces
$J = A$.  If $u$ is the open tripotent corresponding to
${\mathfrak d}$ in $A''$, then $A''(u) = A''$, and so
$u$ is unitary. Since $1 \in {\mathfrak d}$ we have $u^* =
u^* 1 \geq 0$, so that $u = 1$.  Thus ${\mathfrak d} = A_+$.
\end{proof}

{\bf Remark.}  \ The spaces $X$ satisfying the last corollary, which
are also maximally ordered (resp.\
complete in the sense defined just above Theorem \ref{maintro2}) are exactly
the unital operator systems.

\medskip

{\bf Example.}
As a sample illustration of how our results
may be applied in concrete situations,
we  show that if $A'$ is the dual of a C*-algebra $A$,
with its usual cone, then there may exist no isometric positive map
from $A'$ into another C*-algebra.   (A later more elementary proof
of this fact was found in \cite{BNW}.)
We prove it in the case that
$A = \ell^\infty_2$.
In this case  the map $j  : (\alpha,\beta) \mapsto \alpha 1 + \beta z$,
where $z(e^{i \theta}) = e^{i \theta}$,
is a unital complete isometry from $\ell^1_2$ into the C*-algebra $B$
of continuous functions on the unit circle, and in fact
the circle is  well known to be the Shilov boundary
of $\ell^1_2$ (see \cite[Example 4.1.9 (1)]{BLM}),
and so $(B,j)$ is the C*-envelope $C^*_e(A')$.    If
there existed an isometric positive map from $A'$ into another
C*-algebra, then this map would be completely isometric,
since it was noticed by Paulsen that $\ell^1_2$ has exactly one
operator space structure (see \cite[Proposition 3.2]{Pis}
for a simple proof of the latter fact).
Equivalently, the usual cone on $A'$ is an
operator space cone.  By the `nonmatricial cone' case of
Corollary  \ref{DDz},
$(B,j)$ is also the ordered ternary envelope of $A'$.
Hence $j((0,1)) = z$ is a positive function on the unit circle,
which is absurd.

We remark in passing that the fact proved in the last paragraph
shows that the main results about unitizations of ordered spaces in
the paper \cite{K1} are not correct as stated, and this led to the
correction \cite{K2}.

\medskip

We now turn to some other interesting conditions that a cone
might satisfy, some of which also ensure that the cone is
`complete' (that is, the embedding $j : X \to
{\mathcal T}^o(X)$ is a complete order embedding).

Any  natural cone ${\mathfrak c}$ in a TRO $Z$ has the property
that any element in $\overline{{\rm  Span}}({\mathfrak c})$ may
be written as $x = c_1 - c_2 + i(c_3 - c_4)$, where $c_i \in {\mathfrak c}$
and $c_1^* c_2 = 0$ and $c_3^* c_4 = 0$.   (We remark that if $u$ is
the tripotent associated with the cone,
then the product of $c_1$ and $c_2$ in $Z(u)$ is
$c_1 u^* c_2 = u u^* c_1 u^* c_2 = u (c_1^\sharp)^* c_2 =
u c_1^* c_2$, which is $0$ iff $c_1^* c_2 = 0$.)  With this in mind,
it is natural to consider operator space cones
${\mathfrak c}$ in an operator space $X$
with the property that any element in $\overline{{\rm  Span}}({\mathfrak c})$ may
be written as $x = c_1 - c_2+ i(c_3 - c_4)$, where $c_i \in {\mathfrak c}$
and $j(c_1)^* j(c_2) = 0$ and $j(c_3)^* j(c_4) = 0$.  Here $j : X \to {\mathcal T}(X)$ is
the Shilov boundary embedding.   We  say that an operator space cone
is {\em orthogonalizing} if it has this property.

\begin{proposition}  \label{orth}  Let $X$ be an operator space
with an orthogonalizing operator space cone ${\mathfrak c}$.
If ${\mathfrak d}$ is the natural cone
of ${\mathcal T}^o(X)$, then
 ${\mathfrak d} \cap \overline{{\rm  Span}}({\mathfrak c}) = {\mathfrak c}$.  If in addition
${\mathfrak c}$ densely spans $X$, then $X$ is completely order embedded
in ${\mathcal T}^o(X)$, and so $X$ is maximally ordered.
\end{proposition} \begin{proof}
To see this, we view $Z = {\mathcal T}^o(X)$ as a
subTRO in a C*-algebra $B$ and $J(Z)$ as a C*-subalgebra  of $B$.
We are also viewing $\overline{{\rm  Span}}({\mathfrak c}) \subset J(Z)$,
so that, if $x \in {\mathfrak d} \cap \overline{{\rm  Span}}({\mathfrak c})$,
we may write as above $x = c_1 - c_2+ i(c_3 - c_4)$,  with  $c_i \in J(Z)_+$.
Since $x = x^*$ in $J(Z)$,
we must have $x =  c_1 - c_2$.  Since (using the product of $J(Z)$)
 $$0 \leq c_2 x c_2 =
- c_2^3 \leq 0 ,$$
we deduce that $c_2 = 0$ and $x = c_1 \in {\mathfrak c}$.
The last assertion is obvious (using also Corollary \ref{fma2}).  \end{proof}

In a similar spirit, we remark that since the span of
a natural cone  in a  TRO is an inner ideal,
it seems of interest to consider
operator space cones on an operator space $X$ such the
the span of the cone is the analogue of
an `inner ideal' in $X$.  More specifically,
we say that an operator space cone on $X$ is {\em inner}
if $J = \overline{{\rm  Span}}({\mathfrak c})$ is a {\em
generalized quasi-$M$-ideal} of $X$.
The  term  `quasi-$M$-ideal' is due to Kaneda \cite{Kan} (and is
a variant of the one-sided $M$-ideals considered e.g.\ in \cite{BZ}).
By a generalized quasi-$M$-ideal, we mean a subspace $J$ of $X$
such that the weak* closure $J^{\perp \perp}$ of $X''$, viewed as a subspace
of the  ternary envelope $W = ({\mathcal T}(X''),j)$ of $X''$, equals
$p j(X'') q$, where $p$ and $q$ are projections on $W$ which are,
respectively, right and left module maps on $W$.   (In the language
of \cite{BZ} for example,  $p$ and $q$ are orthogonal projections in
${\mathcal A}_\ell(W)$ and ${\mathcal A}_r(W)$ respectively.)
The  generalized quasi-$M$-ideals in a TRO are exactly the
inner ideals (by
e.g.\ the proof in the discussion before Proposition 5.2 in \cite{BHN}).

\begin{theorem}  \label{inn}  Let $X$ be an operator space
with an operator space cone ${\mathfrak c}$ which is inner.
\begin{itemize}  \item [1)]  If ${\mathfrak d}$ is the natural cone
of ${\mathcal T}^o(X)$, then
${\rm Span}({\mathfrak d}) \cap X = \overline{{\rm  Span}}({\mathfrak c})$,
and ${\mathfrak d} \cap X = {\mathfrak d} \cap \overline{{\rm  Span}}({\mathfrak c})$.
 \item [2)]  If also  ${\mathfrak c}$ is orthogonalizing then
${\mathfrak c}$ is complete
(that is, the embedding of $X$ in its ordered ternary envelope ${\mathcal T}^o(X)$
 is a complete order embedding).
\end{itemize} \end{theorem}

\begin{proof}   Write $i_X : X \to X''$ for the canonical injection.
 We use the notation above, so that
$W = ({\mathcal T}(X''),j)$ is the ternary envelope of $X''$.
By \cite[Lemma 5.3]{BHN}, we have that $(\langle j(i_X(X)) \rangle ,
j \circ i_X)$ is a  ternary envelope of $X$.
Here $\langle j(i_X(X)) \rangle$ is the subTRO of $W$ generated by $j(i_X(X))$.
Let  $E = W''$, and $Z = \langle j(i_X(X)) \rangle$, then
$Z'' \cong Z^{\perp \perp}$ is a subWTRO of $E$.
Since $Z$ is the ternary envelope of $X$, there is a natural cone
${\mathfrak d}$ on $Z$ making $Z$ the ordered ternary envelope.
Let $u$ be the associated tripotent in $Z'' \subset E$.
For any $x \in {\mathfrak c}$ we have $p j(i_X(x)) q = j(i_X(x))$.
It follows that $p r(x) q = r(x)$, and so $p u q = u$.
It follows that $p z q = z$ for any $z \in Z(u)$.
Thus $$Z(u) \cap j(i_X(X)) \subset p j(X'') q \cap j(i_X(X)) = j(i_X(J)) .$$
That is, $Z(u) \cap X = J$, and ${\mathfrak d}  \cap X = {\mathfrak d}  \cap J$.

The second part follows from the first part and Proposition \ref{orth}.
\end{proof}

{\bf Remark.}  It seems possible that a converse may hold,
that is,  if the cone is complete then it is inner.

\medskip

{\bf Closing remark.} In the sequel paper \cite{BNW}, we study the
case of operator spaces $X$ which have an involution $*$ and  matrix
cones ${\mathfrak c}_n \subset M_n(X)_{sa}$.  The morphisms in this
category are all $*$-linear, of course. In this case, the ordered
ternary envelope becomes a $*$-TRO, and one must use the ordered
$*$-TRO theory developed in \cite{BW} in place of the ordered TRO
theory in the present paper.   Thus all tripotents $u$ occurring are
also selfadjoint, and central in the sense that $u z = z u$ for all
$z \in Z$.  If  $X$ is an ordered operator space, and if ${\mathcal
T}(X)$ is its `ternary $*$-envelope', then there exists a natural
(in the sense of \cite{BW}) cone on ${\mathcal T}(X)$ containing the
(image of the) cone of $X$, and one can then take the intersection
of all such natural cones to obtain the {\em ordered ternary
$*$-envelope} ${\mathcal T}^o(X)$. The statement of the universal
property of this envelope is similar to that of  Theorem
\ref{maintro2}, but curiously, the proof seems to be completely
different, for the reason that range tripotents need not be central.
Most of our other results from this section, have obvious analogues
in this `selfadjoint case', which we shall not take the time to
spell out.

\medskip

{\bf Acknowledgements.}   We thank Upasana Kashyap for comments on a draft of our
paper.   She has also pointed out to us the validity of the
analogue of \cite[Theorem 4.20]{BW} for
TROs and natural cones in the sense of the present paper.

\end{document}